\documentclass[leqno,11pt]{article}
\usepackage{latexsym}
\usepackage{amssymb}
\usepackage{amsfonts}
\usepackage{amsmath}
\usepackage{color}
\usepackage{epsfig}

\usepackage{caption}
\makeatletter
\long\def\unmarkedfootnote#1{{\long\def\@makefntext##1{##1}\footnotetext{#1}}}
\makeatother

\newtheorem{definition}{Definition}[section]

\newtheorem{lemma}[definition]{Lemma}

\newtheorem{theorem}[definition]{Theorem}

\newtheorem{remark}[definition]{Remark}

\def\m2{|\Omega | /2}
\def\M2{\frac{|\Omega |}{2}}
\def\u+{u_+^*}

\def\-p{\overline{p}}

\def\w0{{W_0^{1,p}(\Omega)}}

\def\R{\mathbb R}
\def\N{\mathbb N}

\newcommand{\med}{{\rm med}}
\def\rn{{{\R}^n}}

\def\B{{\mathbb B}^n}
\def\pimez{\tfrac\pi2}

\newcommand{\hh}{{\cal H}^{n-1}}

\newcommand{\medint}{-\kern  -,395cm\int}
\newcommand{\medintinrigo}{-\kern  -,315cm\int}
\newcommand{\medelle}{-\kern  -,235cm L}
\newcommand{\medellenrigo}{-\kern  -,180cm L}
\newcommand{\qed}{\thinspace\null\nobreak\hfill
\hbox{\vbox{\kern-.2pt\hrule height.2pt
depth.2pt\kern-.2pt\kern-.2pt \hbox to1.8mm {\kern-.2pt\vrule
width.4pt \kern-.2pt\raise1.8mm\vbox to.2pt{} \lower0pt\vtop
to.2pt{}\hfil\kern-.2pt \vrule
width.4pt\kern-.2pt}\kern-.2pt\kern-.2pt \hrule height.2pt
depth.2pt \kern-.2pt}}\par\medbreak}

\title{Poincar\'e trace inequalities in $BV(\mathbb B^n)$ with nonstandard normalization} 
\numberwithin{equation}{section}

\author{
  Andrea Cianchi\\
 {\it Dipartimento di Matematica e Informatica \lq\lq U. Dini", Universit\`a di Firenze}\\ {\it Piazza Ghiberti
27, 50122 Firenze, Italy}
\bigskip
\\
Vincenzo Ferone \\
{\it Dipartimento di Fisica ``E. Pancini"}
\\{\it Universit\`a di Napoli ``Federico II"}
\\ {\it Complesso Monte S.Angelo, Via Cintia, 80126 Napoli, Italy}
\bigskip
\\
Carlo Nitsch \\
{\it Dipartimento di Matematica e Applicazioni ``R. Caccioppoli"}
\\{\it Universit\`a di Napoli ``Federico II"}
\\ {\it Complesso Monte S.Angelo, Via Cintia, 80126 Napoli, Italy}
 \bigskip
\\
Cristina Trombetti \\
{\it Dipartimento di Matematica e Applicazioni ``R. Caccioppoli"}
\\{\it Universit\`a di Napoli ``Federico II"}
\\ {\it Complesso Monte S.Angelo, Via Cintia, 80126 Napoli, Italy}
}

\date{}

\begin{document}
\maketitle
\begin{abstract}\noindent
Extremal functions are exhibited in Poincar\'e  trace inequalities
for functions of bounded variation in the  unit ball $\B$ of the
$n$-dimensional Euclidean space $\rn$. Trial functions are subject
to   either a vanishing mean value condition, or a vanishing median
condition in the whole of $\B$, instead of just on $\partial \B$, as
customary. The  extremals in question take a different form,
depending on  the constraint imposed. In particular, under the
latter  constraint,  unusually shaped extremal functions appear. A
key step in our approach is a characterization of the sharp constant
in the relevant  trace inequalities in any admissible domain $\Omega
\subset \rn$, in terms of an isoperimetric inequality for subsets of
$\Omega$.
\end{abstract}

\unmarkedfootnote {
}

\section{Introduction}\label{intro}

%{\color{blue}}

Let $\Omega$ be a bounded connected open set -- briefly, a domain --
in $\rn$, $n \geq 2$. Under suitable regularity assumptions on its
boundary $\partial \Omega$, a linear trace operator is defined on
the space $BV(\Omega )$ of real-valued functions of bounded
variation in $\Omega$. Such an operator maps any function $u \in
BV(\Omega)$ into a function  $\widetilde u \in L^1(\partial \Omega
)$,  the Lebesgue space of integrable functions on $\partial \Omega$
with respect to the $(n-1)$-dimensional Hausdorff measure $\hh$.
Moreover, the operator
\begin{equation}\label{traceop}
BV(\Omega) \ni u \mapsto \widetilde u \in L^1(\partial \Omega )
\end{equation}
is bounded; namely, there exists a constant $C$ such that
\begin{equation}\label{traceineq}
\| \widetilde u \|_{L^1(\partial \Omega )} \leq C \|u\|_{BV(\Omega)}
\end{equation}
for every $u \in  BV(\Omega)$. Here,
\begin{equation}\label{bvnorm}
 \|u\|_{BV(\Omega)} = \|u\|_{L^1(\Omega)} +  \|Du\|(\Omega),
 \end{equation}
the standard norm in $BV(\Omega)$, where $\|Du\|(\Omega)$ denotes
the total variation in $\Omega$ of the vector-valued Radon measure
$Du$.
\par
Poincar\'e type inequalities also hold, with the full norm
$\|u\|_{BV(\Omega)}$ replaced by just $\|Du\|(\Omega)$, provided
that trial functions $u$ are  normalized by an appropriate operator
$BV(\Omega) \ni u \mapsto T(u) \in \mathbb R$. The relevant
inequalities take the form
\begin{equation}\label{poincineq}
\| \widetilde u - T(u)\|_{L^1(\partial \Omega )} \leq C
\|Du\|(\Omega)
\end{equation}
for  $u \in  BV(\Omega)$. A classical choice for $T(u)$ is the mean
value ${\rm mv}_{\partial \Omega}(\widetilde u)$ of $\widetilde u$
over $\partial \Omega$, given by ${\rm mv}_{\partial
\Omega}(\widetilde u) = \frac 1{\hh (\partial \Omega)} \int
_{\partial \Omega} \widetilde u(x) \,d \hh (x)$. Another customary
option is to take $T(u) = {\rm med}_{\partial \Omega}(\widetilde
u)$, the median of $\widetilde u$ on $\partial \Omega$, defined as
${\med}_{\partial \Omega } (\widetilde u) = \inf \{t \in \mathbb R:
\hh(\{ \widetilde u > t\})\leq \hh (\partial  \Omega)/2\}$.
 General assumptions on the operator $T$ for \eqref{poincineq} to
 hold can be exhibited via a specialization of an abstract result from \cite[Lemma 4.1.3]{Z}.
 \par
 In the present paper we focus on the unconventional cases when either
 \begin{equation*}%\label{mvOmega}
 T(u)=  {\rm mv}_\Omega (u),
 % \quad \hbox{where} \quad {\rm mv}_\Omega (u) = \frac 1{|\Omega|} \int
%_{\Omega}  u(x) \,dx\,
\end{equation*}
where ${\rm mv}_\Omega (u) = \frac 1{|\Omega|} \int _{\Omega}  u(x)
\,dx$, the mean value of $u$ in $\Omega$, or
\begin{equation*}%\label{medOmega}
T(u)=  {\rm med}_\Omega (u),
 %\quad \hbox{where} \quad
%{\med}_{\Omega } ( u) = \inf \{t \in \mathbb R:
%|\{  u > t\}|\leq |\Omega|/2\},
\end{equation*}
where  ${\med}_{\Omega } ( u) = \inf \{t \in \mathbb R: |\{  u
> t\}|\leq |\Omega|/2\}$, the median of $u$ in $\Omega$.  Here,
$|\cdot|$ denotes  Lebesgue measure.
 These choices
make   inequality \eqref{poincineq} nonstandard, in that its
left-hand side  combines quantities depending both on $\widetilde u$
and
 on $u$.
 \par
  We are concerned with the problem of the optimal constant $C$ in
  \eqref{poincineq} for these two choices of $T$.  Our contribution amounts to a characterization of such optimal
 constants  in terms of
 geometric constants of isoperimetric type, and, primarily, to
 the explicit description of the extremal functions in the case when
 $\Omega$ is a ball. In fact, due to the scaling invariance of the relevant inequalities, we shall deal, without loss of generality, with
     the unit ball $\mathbb B^n$, centered at $0$, in $\mathbb R^n$.
    \par Interestingly, the extremals in question differ substantially
 according to whether  $T(u)={\rm
 mv}_{\mathbb B^n} (u)$ or $T(u)={\med}_{\mathbb B^n} ( u)$. As in all known
 sharp  Poincar\'e type inequalities for functions of bounded
 variation, in both cases
 extremal functions are characteristic functions of subsets $\mathbb B^n$.
 %However, %under ${\rm
% mv}_{\mathbb B^n}$ normalization,
As far as the  inequality
\begin{equation}\label{intromv}
\| \widetilde u - {\rm mv}_{\mathbb B^n}(u)\|_{L^1(\partial \mathbb
B^n )} \leq C \|Du\|(\mathbb B^n),
\end{equation}
with an optimal constant $C$, is concerned, extremals turn out to be
 characteristic functions of
 half-balls. Hence,   inequality  \eqref{intromv} shares the
 same  extremals with  the more standard  Poincar\'e  trace
 inequality in $BV(\B)$ with ${\rm
 mv}_{\partial \mathbb B^n}$ normalization (at least for $n \geq 3$)
 \cite{cianchitrace}, and with the mean value
 Poincar\'e inequality in the whole of $\mathbb B^n$
 \cite{Cpoincare}.
 \\ By contrast,
 characteristic functions of   a striking kind of sets
  are extremals
in the inequality
\begin{equation}\label{intromed}
\| \widetilde u - {\rm med}_{\mathbb B^n}(u)\|_{L^1(\partial \mathbb
B^n )} \leq C \|Du\|(\mathbb B^n)
\end{equation}
with optimal constant $C$.
 Unlike the
 case when $T(u) = {\med}_{\partial \mathbb B^n} ( u)$, where
  characteristic functions of half-balls  are still extremals
 \cite{BM}, the extremals in \eqref{intromed} are characteristic functions of half-moon
 shaped subsets of $\mathbb B^n$. In particular, such extremals are not even convex. The isoperimetric nature of
 the optimal constant in  inequality \eqref{intromed}, to
 which we alluded above, helps in accounting for this seemingly
 surprising conclusion.
 \par
  The geometric characterizations of the sharp constant in inequalities \eqref{intromv} and \eqref{intromed}
 are established in Section \ref{geometric}, where some
 definitions and
 basic facts from the theory of functions of bounded variation and
 of sets of finite perimeter are also recalled.
 The results of Section \ref{geometric} are then exploited in Sections \ref{sec2} and \ref{sec1} as  point of departure for the proof of
   sharp forms of inequalities  \eqref{intromv} and \eqref{intromed}, respectively.
Variant  inequalities, where $|\Omega|/2$
 is replaced with an arbitrary fraction of $|\Omega|$ in
 the definition of median, or where trial functions $u$ are required to vanish on a
 subset of $\Omega$ of prescribed measure, are also considered in Sections \ref{geometric} and \ref{sec1}.
\par
We conclude this section by mentioning that   trace inequalities in
Sobolev type spaces, involving optimal constants, have been
extensively investigated in the literature. Contributions along this
line of research include \cite{AFV, AMR, BGP, Br, BrF, Cimosertrace,
CFNT,  DDM, Es1, MV, MV1, Ma0, Matesi, Mazbook, Na, Ro, W}. Sharp
forms of
 Poincar\'e type inequalities  for Sobolev functions and functions of bounded variation,   involving norms
 of $u$ in the whole of $\Omega$,  are the object of
 \cite{BK, BoV, BrV, Cpoincare, DG, DN, EFKNT, ENT, FNT, GW, Le}.
 In particular, the paper \cite{naz} deals with Poincar\'e type
 inequalities for norms  on $\Omega$ of  functions from the Sobolev space $W^{1,2}(\Omega)$, subject to a
 normalization on their traces,
 which are, in a sense, complementary  to those considered here.

\section{Isoperimetric constants}\label{geometric}

Let $E$ be a measurable set in $\rn$. The essential boundary
$\partial ^M E$ of $E$ is defined as the complement in $\mathbb R^n$
of the sets of points of densities $0$ and $1$ with respect to $E$.
One has that $\partial ^M E$  is a Borel set, and $\partial ^M E
\subset \partial E$, the topological boundary of $E$.
\\
The set $E$ is said to be of finite perimeter
 relative to  an open set $\Omega \subset \rn$ if $D\chi _E $, the distributional derivative of the characteristic function $\chi_E$ of $E$, is a
vector-valued Radon measure in $\Omega$ with finite total variation
in $\Omega$.
 The perimeter of $E$  relative to $\Omega$ is defined as
\begin{equation}\label{gen11}
P(E; \Omega ) = \|D \chi _E
\|(\Omega ).
\end{equation}
 A result from geometric measure theory tells us that
$E$ is of finite perimeter in $\Omega$ if and only if $\hh (\partial
^M E \cap \Omega)< \infty$; moreover,
\begin{equation}\label{boundaries3}
P(E; \Omega ) = \hh (\partial ^M E \cap \Omega )
\end{equation}
\cite[Theorem 4.5.11]{Fed}. A domain  $\Omega$ in $\mathbb R^n$ will
be called
  admissible  if
 $\hh (\partial \Omega
 )<\infty $,
 $\hh (\partial \Omega
 \setminus \partial ^M\Omega)=0$, and
 \begin{equation}\label{min}
\min\{\hh (\partial ^M E \cap
\partial \Omega) \, , \hh ( \partial \Omega \setminus \partial ^M E)\} \leq C \hh
(\partial ^M E \cap \Omega)
\end{equation}
for some positive constant $C$ and every measurable set $E \subset
\Omega$ \cite[Definition 5.10.1]{Z}. In particular, any Lipschitz
domain is an admissible domain.
\\
If $\Omega$ is an admissible domain, the boundary trace $\widetilde
u $ of a function  $u \in BV(\Omega )$ is well defined for
$\hh$-a.e. $x \in \partial \Omega$ as
 \begin{equation}\label{trace0}
\widetilde u (x)= \lim _{r \to 0} \frac 1{|B_r(x) \cap
\Omega|} \int _{B_r(x) \cap \Omega}u(y)\, dy\,,
\end{equation}
 where $B_r(x)$ denotes the ball centered at $x$,
with radius $r$  \cite[Corollary 9.6.5]{Mazbook}. The assumption
that $\Omega$ be an admissible domain is necessary and sufficient
for $\widetilde u$ to belong to $L^1(\partial \Omega )$ for every
function $u \in BV(\Omega)$ -- see \cite{AG} and \cite[Theorem
9.5.2]{Mazbook}. Moreover,
 $L^1(\partial \Omega )$ cannot be replaced with any smaller Lebesgue
 space independent of $u$.
 \\ Alternate notions of
 the boundary trace of a function of bounded variation can be found in the literature.
 One definition
  relies upon the notion of upper and lower approximate limits of the
extension of $u$ by $0$ outside $\Omega$ \cite[Definition
5.10.5]{Z}. Another possible definition is that of rough trace in
the sense of \cite[Section 9.5.1]{Mazbook}. Both of them
 coincide with $\widetilde u$, up
to subsets of $\partial \Omega$ of $\hh$-measure zero. \\ If
$\Omega$ is a Lipschitz domain, and the function $u$ enjoys some
additional regularity property, such as membership  to the Sobolev
space $W^{1,1}(\Omega )$, then the trace of $u$ on $\partial \Omega$
 defined as the limit of the restrictions to $\partial
\Omega$ of approximating sequences of smooth functions on $\overline
\Omega$ also  agrees with  $\widetilde u$ for $\hh$-a.e. point on
$\partial \Omega$.

\medskip
\par
Given an admissible domain $\Omega$, let us denote by $C_{\rm mv}(\Omega)$ the optimal constant in the inequality
\begin{equation}\label{tracemv1}
\|\widetilde u - {\rm mv}_\Omega (u)\|_{L^1(\partial \Omega )} \leq C_{\rm mv}(\Omega) \|Du\|(\Omega)
\end{equation}
for  $u \in BV(\Omega)$. Our first result asserts that $C_{\rm mv}(\Omega)$ agrees with the isoperimetric constant
\begin{equation}\label{Ktracemv}
K_{\rm mv}(\Omega) = \sup_{ E\subset\Omega } \frac{|E| \,\hh (\partial \Omega \setminus \partial ^M
E) + |\Omega \setminus E|\,\hh (\partial ^M E \cap
\partial \Omega)}{|\Omega|\, \hh (\partial ^M E \cap \Omega)}.
\end{equation}
Here, and in similar occurrences in what follows, we tacitly assume
that the supremum is extended over non-negligible subsets $E$ of
$\Omega$.

\begin{theorem}\label{Ctracemv}
Let $\Omega$ be an admissible domain in $\rn$, with $n \geq 2$. Then
\begin{equation}\label{gen1}
C_{\rm mv}(\Omega) = K_{\rm mv}(\Omega).
\end{equation}
Equality holds in \eqref{tracemv1} for some nonconstant function $u$
if and only if the supremum is attained in \eqref{Ktracemv} for some
set $E$. In particular, if $E$ is an extremal set in
\eqref{Ktracemv}, then the function $a \chi _E + b$ is an extremal
function in \eqref{tracemv1} for every $a\in \mathbb R\setminus \{
0\}$ and $b \in \R$.
\end{theorem}

\smallskip

\begin{remark}\label{sobolevmv}
{\rm Inequality  \eqref{tracemv1}, with $C_{\rm mv}(\Omega) = K_{\rm
mv}(\Omega)$, holds in particular, for every function $u$ in the
Sobolev space $W^{1,1}(\Omega)$, since the latter is contained in
$BV(\Omega)$. For any such function $u$, the total variation
$\|Du\|(\Omega)$ agrees  with $\|\nabla u\|_{L^1(\Omega )}$, where
$\nabla u$ denotes the weak gradient of $u$. The constant $C_{\rm
mv}(\Omega)$ continues to be optimal in $W^{1,1}(\Omega)$, since any
function $u \in BV(\Omega)$ can be approximated by a sequence of
functions $u_k \in W^{1,1}(\Omega)$ in such a way that
$$\hbox{$\widetilde {u_k} = \widetilde u$ \quad and \quad $\lim _{k \to \infty}\|\nabla u_k\|_{L^1(\Omega )} = \|Du\|(\Omega
)$}.$$  The existence of the sequence $\{u_k\}$  follows, for
instance, from \cite[Theorem 1.17 and Remark 1.18]{Gi}. Of course,
the last part of the statement of Theorem \ref{Ctracemv} does not
apply when dealing with Sobolev functions, since characteristic
functions of subsets of $\Omega$ are not weakly differentiable.
%\\
%In particular, if $\Omega$ is a Lipschitz domain, then the traces of
%functions in $W^{1,1}(\Omega)$ can be classically defined by
%approximation via functions in $C^1(\overline \Omega )$. Hence,
%Theorems \ref{estimate1} and \ref{estimatemean} continue to hold, in
%the class of Lipschitz domains, with $C_{\rm med}(\Omega)$ and
%$C_{\rm mv}(\Omega)$ regarded  as the optimal constants in the trace
%inequalities \eqref{trace} and \eqref{tracemean} for $u \in
%C^1(\overline \Omega )$
}
\end{remark}

\medskip
\par\noindent {\bf Proof of Theorem \ref{Ctracemv}}. Let us begin by showing that
\begin{equation}\label{C<Kmv}
C_{\rm mv}(\Omega) \leq K_{\rm mv}(\Omega),
\end{equation}
namely that
\begin{equation}\label{gen2}
\|\widetilde u - {\rm mv}_\Omega (u)\|_{L^1(\partial \Omega )} \leq C_{\rm mv}(\Omega) \|Du\|(\Omega)
\end{equation}
for  $u \in BV(\Omega)$. Define $u_+ = \frac{u+|u|}2$ and $ u_- = \frac{|u|-u}2$, the positive
and the negative parts of $u$. Since $u = u_+ - u_-$ and $\widetilde u
%= \widetilde{u_+ - u_-}
 = \widetilde {\, u_+ \,} - \widetilde
{\,u_-\,}$,
\begin{equation}\label{gen3}\|\widetilde u - {\rm mv}_\Omega (u)\|_{L^1(\partial \Omega)}\leq \|\widetilde
{\,u_+\,} - {\rm mv}_\Omega (u_+)\|_{L^1(\partial \Omega)} + \|\widetilde {\,u_-\,} -
{\rm mv}_\Omega (u_-)\|_{L^1(\partial
\Omega)}.
\end{equation}
 Moreover, by \cite[Corollary 9.1.2]{Mazbook},
\begin{equation}\label{gen4}
\|Du\|(\Omega) = \|D(u_+)\|(\Omega)+ \|D(u_-)\|(\Omega).
\end{equation}
\par\noindent Thus, we may assume that $u \geq 0$ in \eqref{C<Kmv}, in which case
\begin{equation}\label{gen5}
\widetilde u (x) = \int _0^\infty \chi _{\{\widetilde u \geq t\}}
(x) \, dt \qquad \hbox{for $\hh$-a.e. $x \in \partial \Omega$,}
\end{equation}
%holds for $\hh$-a.e. $x \in \partial \Omega$.
and
\begin{align}\label{gen6}
{\rm mv}_\Omega (u) & = \frac 1{|\Omega|} \int _{ \Omega}u(x)\, dx = \frac 1{|\Omega|} \int _{ \Omega}\bigg(\int
_0^\infty \chi
_{\{ u \geq t\}} (x) \, dt\bigg)\, d x \\
\nonumber &= \frac 1{|\Omega|} \int _0^\infty
\bigg(\int _{ \Omega}\chi _{\{ u \geq t\}} (x)\,
d x\bigg) \, dt = \frac 1{|\Omega|}  \int
_0^\infty |\{u \geq t\}|\, dt.
\end{align}
Hence, the following chain holds:
\begin{align}\label{tracemv1bis}
\|\widetilde u -  {\rm mv}_\Omega (u)\|_{L^1(\partial \Omega)} & = \int
_{\partial \Omega} \bigg|\widetilde u (x) - \frac 1{|\Omega|} \int
_\Omega u(y)\, dy|\bigg| d\hh (x)
\\ \nonumber & = \int _{\partial \Omega}
\bigg|\int _0^\infty \chi_{\{\widetilde u \geq t\}}(x)\, dt - \frac
1{|\Omega|} \int_0^\infty |\{u \geq t\}|dt\bigg| d\hh (x)
\\ \nonumber & \leq
\int _{\partial \Omega} \int _0^\infty \bigg| \chi_{\{\widetilde u
\geq t\}}(x) - \frac 1{|\Omega|} \int_0^\infty |\{u \geq t\}|\bigg|
dt\, d\hh (x)
 \\ \nonumber & = \int _0^\infty \int _{\partial \Omega} \bigg|
\chi_{\{\widetilde u \geq t\}}(x) - \frac 1{|\Omega|} \int_0^\infty
|\{u \geq t\}|\bigg| d\hh (x)\, dt
\\ \nonumber  & =
\frac 1{|\Omega|} \int _0^\infty \big[\hh (\{\widetilde u \geq t\} )
(|\Omega | - |\{u\geq t\}|)%\, dt
\\ \nonumber  & \quad \quad \quad+
%\frac 1{|\Omega|} \int _0^\infty
(\hh (\partial \Omega ) - \hh (\{\widetilde u \geq t\}))|\{u\geq
t\}|\big]\, dt.
\end{align}
One has that,
\begin{equation}\label{trace8}
\hh (\{\widetilde u \geq t\}) = \hh (\partial ^M \{ u \geq t\} \cap
\partial \Omega) \quad \hbox{for a.e. $t>0$}
\end{equation}
(see e.g. \cite[Equation (2.6)]{cianchitrace}). On the other hand,
by  \eqref{Ktracemv},
\begin{multline}\label{gen7}
\hh (\partial ^M \{ u \geq t\} \cap \partial \Omega)
(|\Omega | - |\{u\geq t\}|) + \hh (\partial \Omega \setminus
\partial ^M \{ u \geq t\})|\{u\geq t\}| \\ \leq   K_{\rm mv}(\Omega) |\Omega|\, \hh (\partial ^M \{ u \geq t\} \cap \Omega) \quad \hbox{for a.e. $t>0$.}
\end{multline}
From \eqref{tracemv1bis}, \eqref{trace8} and \eqref{gen7} we infer that
\begin{align}\label{gen8}
\|\widetilde u -  {\rm mv}_\Omega (u)\|_{L^1(\partial \Omega)}  & =
\frac 1{|\Omega|} \int _0^\infty \big[\hh (\partial ^M \{ u \geq t\}
\cap \partial \Omega) (|\Omega | - |\{u\geq t\}|)
\\ \nonumber & \quad \quad \quad + (\hh (\partial
\Omega ) - \hh (\partial ^M \{ u \geq t\} \cap \partial
\Omega))|\{u\geq t\}|\big]\, dt
\\ \nonumber  & = \frac 1{|\Omega|} \int _0^\infty \big[\hh (\partial ^M \{ u \geq t\} \cap \partial \Omega)
(|\Omega | - |\{u\geq t\}|)%\,dt
\\ \nonumber  & \quad \quad \quad+ %\frac 1{|\Omega|} \int _0^\infty
\hh (\partial \Omega \setminus
\partial ^M \{ u \geq t\})|\{u\geq t\}|\big]\,
dt
\\ \nonumber  & \leq  K_{\rm mv}(\Omega) \int _0^\infty \hh (\partial ^M \{ u \geq t\} \cap \Omega)\,
dt.
 \end{align}
Finally, the coarea formula
 for $BV$-functions \cite[Theorem 5.4.4]{Z} tells us that
\begin{equation}\label{trace11}
\int _0^\infty \hh (\partial ^M \{ u \geq
 t\}\cap \Omega) \,dt = \|Du\|(\Omega).
 \end{equation}
Combining equations \eqref{gen8} and \eqref{trace11} yields inequality \eqref{C<Kmv}.
\\
In order to prove the reverse inequality in \eqref{C<Kmv}, namely
that
\begin{equation}\label{C>Kmv}
C_{\rm mv}(\Omega) \geq K_{\rm mv}(\Omega),
\end{equation}
consider any set $E \subset \Omega$ of finite perimeter in $\Omega$.
Since, by \eqref{trace8}, $\widetilde{\chi _E} = \chi_{\partial ^M E
\cap \partial \Omega}$ outside a set of $\hh$ measure zero on
$\partial \Omega$,
 \begin{align}\label{gen10}
\|\widetilde{\chi _E} - {\rm mv}_\Omega (\chi _E)\|_{L^1(\partial
\Omega )} & = \int _{\partial \Omega}\Big|\chi_{\partial ^M E \cap
\partial \Omega} (x) - \tfrac{|E|}{|\Omega|}\Big|d \hh (x)
\\ \nonumber
& = \frac 1{|\Omega|}\big(|E| \,\hh (\partial \Omega \setminus
\partial ^M E) + |\Omega \setminus E|\,\hh (\partial ^M E \cap
\partial \Omega)\big).
\end{align}
On the other hand, by \eqref{gen11} and \eqref{boundaries3},
\begin{equation}\label{gen12}
\|D\chi _E\|(\Omega) = \hh (\partial ^M E \cap \Omega).
\end{equation}
Hence, the choice of trial functions $u$ of the form $\chi _E$ in
inequality \eqref{tracemv1} tells us that
\begin{align}\label{gen13}
|E| \,\hh (\partial \Omega \setminus \partial ^M
E) + |\Omega \setminus E|\,\hh (\partial ^M E \cap
\partial \Omega) \leq C_{\rm mv}(\Omega) |\Omega| \hh (\partial ^M E \cap \Omega)
\end{align}
for every set $E$ of finite perimeter in $\Omega$, whence \eqref{C>Kmv} follows.
\\ Assume now that $E$ is any  set at which the supremum is attained in \eqref{Ktracemv}. Thus,
\begin{align}\label{gen14}
K_{\rm mv}(\Omega) & = C_{\rm mv}(\Omega) \geq \frac{\|\widetilde{\chi _E} -
{\rm mv}_\Omega (\chi _E)\|_{L^1(\partial \Omega )}}{ \|D\chi _E\|(\Omega)}
\\ \nonumber & =  \frac{|E| \,\hh (\partial \Omega \setminus \partial ^M
E) + |\Omega \setminus E|\,\hh (\partial ^M E \cap
\partial \Omega)}{|\Omega|\, \hh (\partial ^M E \cap \Omega)} = K_{\rm mv}(\Omega).
\end{align}
Consequently, equality holds in the inequality in \eqref{gen14}.
This means that $\chi _E$, and hence $a\chi _E +b$ for every $a \neq
0$ and $b \in \mathbb R$, is an extremal in \eqref{tracemv1}.
Conversely, assume that $u$ is an extremal in \eqref{tracemv1}, i.e.
$u$ is nonconstant, and equality holds in \eqref{tracemv1}. A close
inspection of the proof of \eqref{C<Kmv} reveals that equality must
hold in the inequality in \eqref{gen8}, applied with $u$ replaced
with $u_+$ and $u_-$. Hence, equality has to hold in \eqref{gen7},
with $u$ replaced with $u_+$ and $u_-$, for a.e. $t \geq 0$. This
tells us that the sets $\{u_{\pm} \geq t\}$ are extremals in
\eqref{Ktracemv} for a.e. $t \in [0, {\rm esssup}\, u_{\pm})$. \qed

\bigskip
\par

We next take into account inequalities with ${\rm med}_{\Omega}$ normalization. Let us call $C_{\rm med}({\Omega})$
 the optimal constant in the inequality
\begin{equation}\label{infinf}
\|\widetilde u - {\rm med}_{\Omega} (u)\|_{L^1(\partial \Omega )} \leq C_{\rm med}({\Omega}) \|Du\|(\Omega)
\end{equation}
for  $u \in BV(\Omega)$. The isoperimetric constant which now comes into play is defined as
\begin{equation}\label{Komegasup}
K_{\rm med}({\Omega}) = \sup_{
\begin{tiny}
 \begin{array}{c}{E\subset\Omega} \\
|E|\le| \Omega|/2
 \end{array}
  \end{tiny}}
\frac{ \hh (\partial ^M E \cap
\partial \Omega) }{\hh (\partial ^M E \cap \Omega)}.
\end{equation}

\begin{theorem}\label{Comega}
Let $\Omega$ be an admissible domain in $\rn$, with $n \geq 2$. Then
\begin{equation}\label{admiss}
C_{\rm med}({\Omega}) = K_{\rm med}({\Omega}).
\end{equation}
Equality holds in \eqref{infinf} for some nonconstant function $u$  if and only if the supremum is attained  in \eqref{Komegasup}
for some set $E$. In particular, if $E$ is an
extremal in \eqref{Komegasup}, then the function $a \chi _E + b$ is an
extremal in \eqref{infinf} for every $a \in \R \setminus \{0\}$ and $b \in \R$.
\end{theorem}

Theorem \ref{Comega} is a special case of a slightly more general
result, which is the content of Theorem \ref{Comegasigma} below. Its
statement involves the following definitions.  Given $\sigma \in (0,
1)$,  set
\begin{equation}\label{medsigma}{\med}_{\Omega, \sigma  } ( u) = \inf \{t \in \mathbb R:
|\{  u > t\}|\leq \sigma |\Omega|\},
\end{equation}
and denote by $C_{\rm med}({\Omega,\sigma})$ the optimal constant in
the trace inequality
\begin{equation}\label{infsigma}
\|\widetilde u - {\rm med}_{\Omega, \sigma} (u)\|_{L^1(\partial \Omega )} \leq C_{\rm med}({\Omega,\sigma}) \|Du\|(\Omega)
\end{equation}
for  $u \in BV(\Omega)$. Moreover,  define
\begin{equation}\label{Komegasigma}
K_{\rm med}({\Omega, \sigma}) = \sup_{
\begin{tiny}
 \begin{array}{c}{E\subset\Omega} \\
|E|\le| \sigma| \Omega|
 \end{array}
  \end{tiny}}
\frac{ \hh (\partial ^M E \cap
\partial \Omega) }{\hh (\partial ^M E \cap \Omega)}.
\end{equation}

\medskip

\begin{theorem}\label{Comegasigma}
Let $\Omega$ be an admissible domain in $\rn$, with $n \geq 2$, and
let $\sigma \in (0,1)$. Set $\rho=\max\{\sigma, 1-\sigma\}$. Then
\begin{equation}\label{admisssigma}
C_{\rm med}({\Omega, \sigma}) = K_{\rm med}({\Omega, \rho}).
\end{equation}
Equality holds in \eqref{infsigma} for some nonconstant function $u$
if and only if the supremum is attained  in \eqref{Komegasigma} for
some set $E$. In particular, if $E$ is an extremal set in
\eqref{Komegasigma}, then the function $a \chi _E + b$ is an
extremal in \eqref{infsigma} for every $a \in \R \setminus \{0\}$
and $b \in \R$.
\end{theorem}
\par\noindent
{\bf Proof}. We begin by proving that
\begin{equation}\label{C<K}
C_{\rm med}({\Omega, \sigma}) \leq K_{\rm med}({\Omega, \rho}).
\end{equation}
Inequality \eqref{C<K} will follow if we show that
\begin{equation}\label{inf0}
\|\widetilde u \|_{L^1(\partial \Omega )}
 \le
 K_{\rm med}(\Omega) \|Du\|(\Omega)
\end{equation}
for every $u \in BV(\Omega )$ such that ${\rm med}_{\Omega, \sigma} (u) =0$. For
any such $u$,
\begin{equation}
|\{x\in \Omega : u_\pm (x) \geq t\}| \leq \rho |\Omega| \quad
\hbox{for $t>0$.}
\end{equation}
%Since $\widetilde u
% = \widetilde {\, u_+ \,} - \widetilde
%{\,u_-\,}$, we have that
Furthermore,
\begin{equation}\label{trace3}\|\widetilde u \|_{L^1(\partial \Omega)}\leq \|\widetilde
{\,u_+\,} \|_{L^1(\partial \Omega)} + \|\widetilde {\,u_-\,}
\|_{L^1(\partial \Omega)}.
\end{equation}
By \eqref{trace3} and \eqref{gen4},
it suffices to prove inequality \eqref{inf0} in the case when $u
\geq 0$ and
\begin{equation}\label{trace5}
|\{x\in \Omega : u (x) \geq t\}| \leq \rho |\Omega| \quad
\hbox{for $t>0$.}
\end{equation}
Let $u$ be any function in $BV(\Omega)$ satisfying these properties.
Owing to \eqref{trace8},
\begin{align}\label{trace7bis}
\|\widetilde u\|_{L^1(\partial \Omega)} & = \int _{\partial \Omega}
\widetilde u (x)\, d\hh (x) = \int _{\partial \Omega} \int _0^\infty
\chi _{\{\widetilde u \geq t\}} (x) \, dt d\hh (x)
\\ \nonumber & =
\int _0^\infty \int _{\partial \Omega}\chi _{\{\widetilde u \geq
t\}} (x)\, \hh (x) \, dt
 = \int _0^\infty \hh(\{\widetilde u \geq
t\})\, dt
\\ \nonumber & = \int _0^\infty \hh
(\partial ^M \{ u \geq t\} \cap
\partial \Omega)\, dt.
\end{align}
On the other hand, by \eqref{trace5} and \eqref{Komegasigma},
\begin{equation}\label{trace9}
\hh (\partial ^M \{ u \geq t\} \cap
\partial \Omega)  \leq K_{\rm med}({\Omega, \rho}) \hh (\partial ^M \{ u \geq t\} \cap
\Omega) \quad \hbox{for a.e. $t>0$.}
\end{equation}
Coupling \eqref{trace7bis} with \eqref{trace9} yields
\begin{equation}\label{trace10}
\|\widetilde u\|_{L^1(\partial \Omega)} \leq K_{\rm med}({\Omega, \rho}) \int
_0^\infty \hh (\partial ^M \{ u \geq t\} \cap \Omega) \, dt.
\end{equation}
 Inequality \eqref{inf0} follows from \eqref{trace10} and
\eqref{trace11}.
\\
We next prove that
\begin{equation}\label{C>K}
C_{\rm med}({\Omega, \sigma}) \geq K_{\rm med}({\Omega, \rho}).
\end{equation}
Let $E \subset \Omega$ be such that $\hh (\partial ^M E \cap \Omega
) < \infty$ and  $|E| \leq \rho|\Omega|$.  Then, either ${\rm
med}_{\Omega, \sigma}(\chi _E) =0$,  or  ${\rm med}_{\Omega,
\sigma}(-\chi _E) =0$, according to whether $\sigma \geq \tfrac 12$
or $\sigma \leq \tfrac 12$. Since $\pm \chi _E \in BV(\Omega)$,
either  $u= \chi _E$ or $u= -\chi _E$ is an admissible trial
function in \eqref{infsigma}. By \eqref{trace7bis},
 \begin{equation}\label{trace12}
\|\widetilde {\pm\chi _E}\|_{L^1(\partial \Omega)}=   \hh (\partial ^M
E \cap
\partial \Omega ).
\end{equation}
From \eqref{infsigma}, \eqref{trace12} and \eqref{gen12} we deduce
that
 \begin{equation}\label{trace14}
\hh (\partial ^M E \cap
\partial \Omega ) \leq C_{\rm med}({\Omega, \sigma}) \hh (\partial ^M E \cap \Omega ),
\end{equation}
whence \eqref{C>K} follows.
\\
Now, assume that $E$ is any set with $\hh (\partial ^M E \cap \Omega
) < \infty$ and  $|E| \leq \rho|\Omega|$, at which equality is
attained in \eqref{Komegasigma}. In particular, either  ${\rm
med}_{\Omega, \sigma}(\chi _E) =0$, or  ${\rm med}_{\Omega,
\sigma}(- \chi _E) =0$. Thus, owing to \eqref{trace12} and
\eqref{gen12},
\begin{align}\label{trace15}
K_{\rm med}({\Omega, \rho}) = C_{\rm med}({\Omega, \sigma}) \geq  \frac{\hh
(\partial ^M E \cap
\partial \Omega )}{\hh (\partial ^M E \cap \Omega )} = K_{\rm med}(\Omega,\rho).
\end{align}
This shows that equality holds in the inequality in \eqref{trace15}.
Therefore, either $\chi _E$, or  $-\chi _E$ is an extremal function
in \eqref{infsigma}, and hence $a\chi _E + b$ is an extremal
function for every $a \in \R \setminus  \{0\}$ and $b \in \R$.
\\
Conversely, assume that equality holds in \eqref{infsigma} for some
function $u \in BV(\Omega)$. We may clearly assume that ${\rm
med}_{\Omega, \sigma} =0$. It is easily seen via an inspection of
the proof of \eqref{infsigma} that then equality must hold in
\eqref{trace9}, with $u$ replaced by $u_+$ and $u_-$, for a.e.  $t
>0$. Hence, the sets $\{u_\pm \geq t\}$ are extremals in
\eqref{Komegasigma} for a.e. $t \in [0, {\rm ess sup} \,u_\pm)$.
\qed

The constant given by
\eqref{Komegasigma} also enters in a trace inequality for
functions subject to a different normalization.

\begin{theorem}\label{Comegazero}
Let $\Omega$ be an admissible domain in $\rn$, with $n \geq 2$, and
let $\sigma \in (0, 1)$.  Then
\begin{equation}\label{zerosigma}
\|\widetilde u \|_{L^1(\partial \Omega )}
 \le
 K_{\rm med}({\Omega, \sigma}) \|Du\|(\Omega)
\end{equation}
for every $u \in BV(\Omega )$ such that
$$|\{u =0\}| \geq \sigma |\Omega|.$$
Equality holds in \eqref{zerosigma} for some function $u$ which does
not vanish identically if and only if equality holds in
\eqref{Komegasigma} for some set $E$. In particular, if $E$ is an
extremal set in \eqref{Komegasigma}, then the function $a \chi _E $
is an extremal in \eqref{zerosigma} for every $a \in \R \setminus
\{0\}$.
\end{theorem}

The proof of Theorem \ref{Comegazero} is completely analogous to (and even simpler than) that of Theorem \ref{Comegasigma}, and will be omitted.

\medskip
\par\noindent

\begin{remark}\label{sobolevmed}
{\rm Considerations as in Remark \ref{sobolevmv}  hold in connection
with Theorems \ref{Comega}, \ref{Comegasigma} and \ref{Comegazero}
as well. These results thus provide a geometric characterization of
the optimal constant in the pertaining trace inequalities also for
functions from the Sobolev space $W^{1,1}(\Omega)$.}
\end{remark}

\section{A trace inequality on $\mathbb B^n$ with mean value normalization}\label{sec2}

In the remaining part of this paper, the geometric characterizations
of the sharp constants in the Poincar\'e trace inequalities provided
by Theorems \ref{Ctracemv} and \ref{Comegasigma} are specialized to
the case when the ground domain $\Omega$ agrees with the ball $\B$.
Its peculiar geometry enables us to exhibit the extremal subsets in
the associated isoperimetric problems, and hence the extremal
functions in the relevant trace inequalities. In the light of
Remarks \ref{sobolevmv} and \ref{sobolevmed}, the resulting
inequalities are not only sharp in $BV (\B)$, but also in
$W^{1,1}(\B)$.
\\
This section is devoted to the problem of the optimal constant
$C_{\rm mv}(\B)$ in the mean value inequality
\begin{equation}\label{Cmvoptn}
\|\widetilde u - {\rm mv}_{\B} (u)\|_{L^1(\partial \B )} \leq C_{\rm
mv}(\B)\|Du\|(\B)
\end{equation}
for   $u \in BV(\B)$. Its solution reads as follows.

\begin{theorem}\label{isopconstpmv} Let $n \geq 2$. Then
%\begin{equation}\label{Cmvoptn}
%\|\widetilde u - {\rm mv}_{\B} (u)\|_{L^1(\partial \B )} \leq
%\frac{n\omega_n}{2\omega_{n-1}}\|Du\|(\B)
%\end{equation}
%for every $u \in BV(\B)$.
\begin{equation*}
C_{\rm mv}(\B)=\frac{n\omega_n}{2\omega_{n-1}}.
\end{equation*}
Equality holds in \eqref{Cmvoptn} if u agrees with the
characteristic function of a half-ball.
\end{theorem}

Theorem \ref{isopconstpmv} is a straightforward consequence of
Theorem \ref{Ctracemv} and of Theorem \ref{isopconstmv} below.

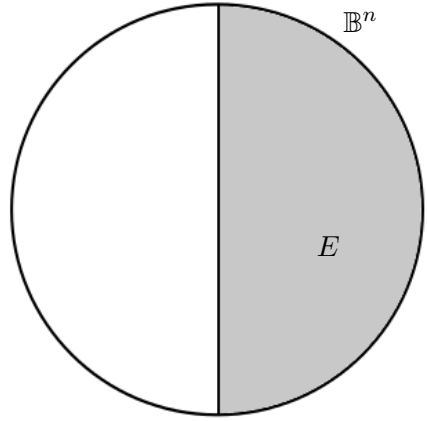
\begin{figure}\setlength\fboxsep{30pt} \setlength\fboxrule{0.0pt}
\hskip3.6cm\fbox{\includegraphics[width=2.2 in]{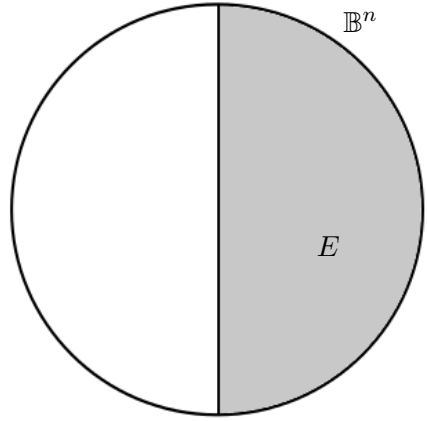}}
\caption{Half-balls are extremal sets for $K_{\rm
mv}(\B)$}\label{intersection_h}
\begin{picture}(0,0)(150,-30)
\put(410,180){$\B$}
\put(400,95){$E$}
\end{picture}
\end{figure}

\begin{theorem}\label{isopconstmv} Let $n \geq 2$. Then
\begin{equation}\label{Kmvoptnew}
K_{\rm mv}(\B)=\frac{n\omega_n}{2\omega_{n-1}}.
\end{equation}
Half-balls are extremal sets for $K_{\rm mv}(\B)$ (see Figure \ref{intersection_h}).
\end{theorem}

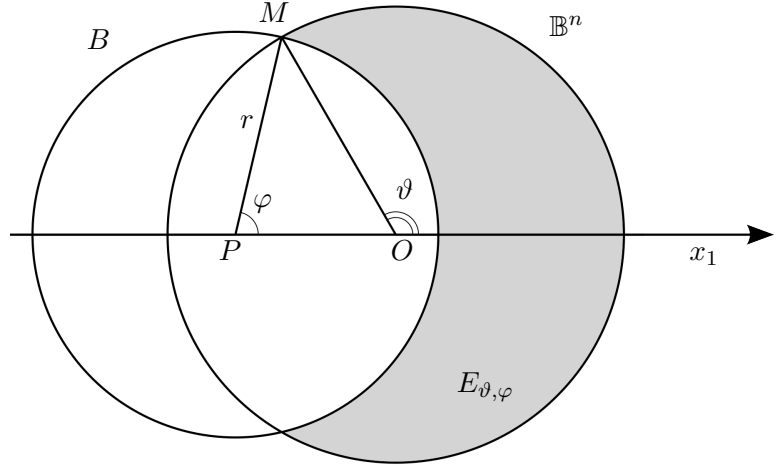
\begin{figure}
\setlength\fboxsep{30pt} \setlength\fboxrule{0.0pt}
\hskip2.5cm\fbox{\includegraphics[width=4 in]{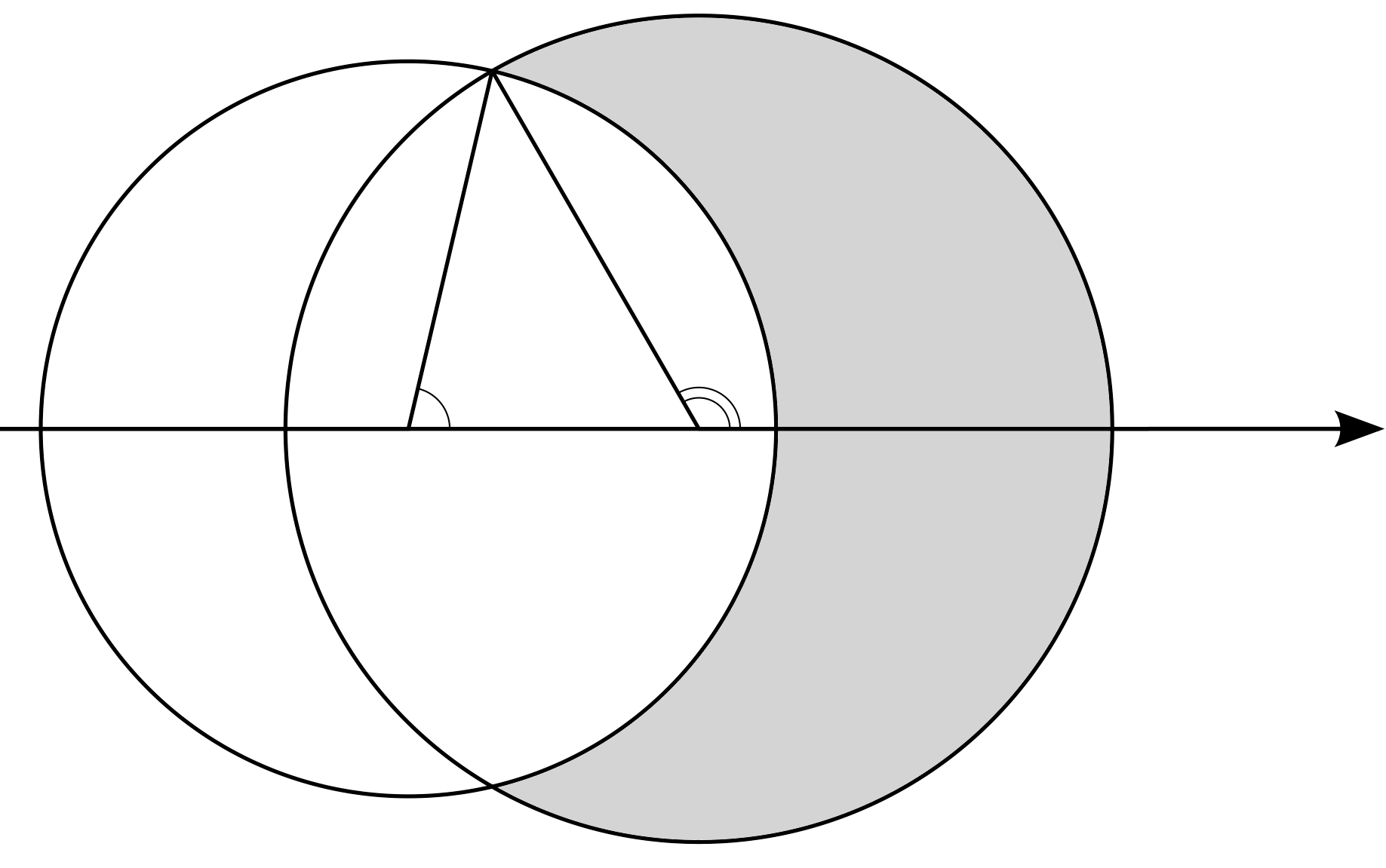}}
\caption{The set $E_{\vartheta,\varphi}$}\label{intersection}
\begin{picture}(0,0)(150,-30)
\put(280,195){$B$} \put(345,205){$M$} \put(455,200){$\B$}
\put(420,65){$E_{\vartheta, \varphi}$} \put(397,138){$\vartheta$}
\put(343,135){$\varphi$} \put(395,115){$O$} \put(330,115){$P$}
\put(508,115){$x_1$} \put(338,165){$r$}
\end{picture}
\end{figure}

Symmetrization, and other ad hoc geometric arguments, enable us to
restrict the analysis of possible extremal sets for $K_{\rm mv}(\B)$
to a two-parameter family of subsets of $\B$, which are the
complement in $\B$ of another ball $B$.
\\
Specifically, let $O$ and $P$ be the centers of $\B$ and $B$,
 and, with reference to Figure \ref{intersection} (for $n=2$), let $E_{\vartheta, \varphi}=\B\setminus
 B$,
where $\vartheta$ denotes the angle between the positive
$x_1$-half-axis and the radius of   $\B$ issued from a point $M \in
\partial \B \cap \partial B$,  and  $\varphi$ denotes the angle between
the same half-axis and the radius of   $B$ through $M$.  The couple
$(\vartheta, \varphi)$ belongs to the set
\begin{equation}\label{dominio}
\Upsilon = \{(\vartheta, \varphi): 0<\vartheta < \pi , 0\leq \varphi
< \vartheta \}.
\end{equation}
The endpoint case when $\varphi =0$ corresponds to the borderline
situation when $B$ is a half-space, and hence $E_{\vartheta, 0}$ is
a spherical segment.
\\
In fact, in the proof of Theorem  \ref{isopconstmv}, we shall only
need to consider sets $E_{\vartheta, \varphi}$ such that $B
\nsubseteq \B$, and
$$\hh(\partial ^M E_{\vartheta, \varphi} \cap\partial \B)\ge
\hh(\partial \B\setminus\partial ^M E_{\vartheta, \varphi}),
$$
namely   couples $(\vartheta,\varphi)$ from the set
\begin{equation}\label{defdelta}
\Theta =\{(\vartheta,\varphi)\in \R^2:\pimez\le \vartheta<\pi
,\>0\le \varphi< \vartheta\}.
\end{equation}
Denote by $r$ the radius of $B$, and observe that, if $\varphi>0$,
then
$$r =  \frac{\sin \vartheta }{\sin \varphi}.$$
Relevant geometric quantities associated with the set $E_{\vartheta,
\varphi}$ can be expressed in terms of the functions $\Psi_{k}$ and
$\Phi_{k}$ defined,
 for $k\in \N\cup\{0\}$, by
\begin{equation}\label{def}
\Psi_{k}(t)=\int_0^{t } \sin ^{k}\tau\,d\tau ,
\end{equation}
and
\begin{equation}\label{defphi}
\Phi_{k}(t)=\int_0^{t } \cos ^{k}\tau\,d\tau ,
\end{equation}
for $t \in [0, \pi]$. The following equations are easily verified:
\begin{equation}\label{ott1}
\Psi_k(t)=\frac{k-1}k\Psi_{k-2}(t)-\frac1k\cos t\,\sin^{k-1}t,
\end{equation}
 \begin{equation}\label{nov24}
 \Phi_k(t)=\frac{k-1}k\Phi_{k-2}(t)+\frac1k\cos^{k-1}t\sin t,
 \end{equation}
for $k\ge2$ and $t \in [0, \pi]$.
%%
%
%
%
%
%
%
%\begin{figure}
%\setlength\fboxsep{30pt} \setlength\fboxrule{0.0pt}
%\hskip1.9cm\fbox{\includegraphics[width=3.5 in]{disegno_cerchi.png}}
%\caption{The set $E_{\vartheta,\varphi}$}\label{intersection}
%\begin{picture}(0,0)(150,-30)
%\put(260,200){$B$} \put(337,215){$M$} \put(455,200){$\B$}
%\put(420,65){$E_{\vartheta, \varphi}$} \put(390,140){$\vartheta$}
%\put(334,137){$\varphi$} \put(385,118){$O$} \put(320,118){$P$}
%\put(490,118){$x_1$} \put(325,165){$r$}
%\end{picture}
%\end{figure}
\\ Computations show that
\begin{align}%\label{eq:1}
\qquad&\hh (\partial ^M E_{\vartheta, \varphi} \cap\partial \B)
=(n-1)\omega_{n-1}\Psi_{n-2}(\vartheta), %\qquad(\vartheta,\varphi)\in\Theta,
&
\label{def1}\\ \notag\\
\qquad&\hh (\partial ^M E_{\vartheta, \varphi} \cap
\B)=(n-1)\omega_{n-1}\Psi_{n-2}(\varphi)\Big(\frac{\sin \vartheta
}{\sin
\varphi}\Big)^{n-1}  , %\qquad(\vartheta,\varphi)\in\Theta,
& \label{def2}
\end{align}
for $(\vartheta,\varphi)\in \Upsilon$, where $\omega_{n}=|\B|$,
namely $\omega_n = \frac{\pi^{n/2}}{\Gamma(1+n/2)}$. Moreover,
\begin{equation}\label{def3}
|E_{\vartheta, \varphi}|=
%\omega_{n-1}(\Psi_{n}(\vartheta)R^{n}-\Psi_{n}(\varphi)r^{n})=
\omega_{n-1}\left(\Psi_{n}(\vartheta)-\Psi_{n}(\varphi)\left(\frac{\sin\vartheta}{\sin\varphi}\right)^{n}\right),
%\qquad(\vartheta,\varphi)\in\Theta,
\end{equation}
and
\begin{equation}\label{defmv3}
|\B \setminus E_{\vartheta ,
\varphi}|=\omega_{n-1}\left(\Psi_{n}(\pi-\vartheta)+\Psi_{n}(\varphi)\left(\frac{\sin\vartheta}{\sin\varphi}\right)^{n}\right).  %\qquad(\vartheta,\varphi)\in\Theta.
\end{equation}
Note that equations \eqref{def1}--\eqref{defmv3} also hold for
$\varphi =0$, which corresponds to the case when $E_{\vartheta , 0}$
is the intersection of $\B$ with a half-space, provided that their
right-hand sides are extended by continuity.
\\
The following equations  will  be used below without further
mentioning:
\begin{equation}\label{nov10}
\omega _n = \omega _{n-1}\Psi _n(\pi ),
\end{equation}
\begin{equation}\label{nov11}
n \omega _n = (n-1) \omega _{n-1}\Psi _{n-2}(\pi ),
\end{equation}
\begin{equation}\label{nov12}
\Psi _k (\pi ) - \Psi_{k} (t ) = \Psi _k (\pi - t),
\end{equation}
\begin{equation}\label{nov13}
\Psi _k (t ) - \Psi_{k} (\pi - \vartheta ) = 2\Phi _k (t - \pimez)
\end{equation}
for $k \geq 0$ and $t \in [0, \pi]$.

Given a set $E \subset \B$, let us  denote by ${\cal Q}_{\rm mv}(E)$
the quotient appearing in \eqref{Ktracemv} for $\Omega = \B$,
 namely
 \begin{equation}\label{defQmv} {\cal Q}_{\rm mv}(E)=\frac{|E| \,\hh
(\partial \B \setminus \partial ^M E) + |\B \setminus E|\,\hh
(\partial ^M E \cap
\partial \B)}{|\B|\, \hh (\partial ^M E \cap \B)}.
\end{equation}
In particular, owing to \eqref{def1}--\eqref{def3},
\begin{equation}
\label{QEmv} {\cal Q}_{\rm mv}(E_{\vartheta, \varphi})
=\frac{\displaystyle n\omega_{n}^2-4(n-1) \omega_{n-1}^2
\Phi_{n-2}(\vartheta-\pimez)\left(\Phi_{n}(\vartheta-\pimez)-\Psi_{n}(\varphi)
\frac{\sin^n\vartheta}{\sin^n\varphi}\right)}{\displaystyle2(n-1)\omega_{n-1}\omega_{n}\Psi_{n-2}(\varphi)\frac{\sin^{n-1}\vartheta}{\sin^{n-1}\varphi}}
\end{equation}
for $(\vartheta, \varphi) \in \Theta$, where the expression on the
right-hand side is extended by continuity for $\varphi =0$.

The following technical lemma will be needed in our proof of Theorem
\ref{isopconstmv}.

\begin{lemma}\label{symmmv1}
Let $n\ge2$, and let
$$\Lambda=\{(t,s)\in \R^2:0\leq  t <
\pimez ,\>0\leq  s< t+\pimez\}.$$ Define  the function $F : \Lambda
\to \mathbb R$ as
\begin{equation}
\label{deffmv} F(t,s)=\displaystyle 1-4\frac{n-1}n
\frac{\omega_{n-1}^2}{\omega_{n}^2}
\Phi_{n-2}(t)\left(\Phi_{n}(t)-\Psi_{n}(s)\frac{\cos^nt}{\sin^ns}\right)-(n-1)\Psi_{n-2}(s)
\frac{\cos^{n-1}t}{\sin^{n-1}s}
\end{equation}
for $(t,s) \in \Lambda$, where the right-hand side is extended by
continuity for $s=0$. Then
\begin{equation}
\label{ineqmv} F(t,s) \leq F(0,0) =0 \quad \hbox{for  $(t,s)\in
\Lambda$ ,}
\end{equation}
and the equality holds in the first inequality only if
$(t,s)=(0,0)$.
%\textcolor[rgb]{0.00,0.00,1.00}{
Hence, $F$ attains its maximum in $\Lambda$ only at  $(0,0)$.
\end{lemma}
{\bf Proof}. Assume first that $t=0$. One has that
$$F(0, s) = 1 - (n-1)
\frac{\Psi_{n-2}(s)}{\sin^{n-1}s}, \qquad\hbox{for $s \in (0,
\pimez]$.}$$ If $g : [0, \pimez] \to \mathbb R$ is the function
defined by $g(s)= \sin^{n-1}s - (n-1) \Psi_{n-2}(s)$ for $s \in [0,
\pimez]$, then $g(0)=0$ and $g'(s) < 0$ for $s \in (0, \pimez]$.
Hence, inequality \eqref{ineqmv} follows for $t=0$.
\\ Next, fix any $t\in(0,\pimez)$.
Equation \eqref{ott1}, with $k=n$, ensures that
\begin{align}\label{nov20}
\frac{\partial F}{\partial s}(t,s) = (n-1) \frac{\cos ^{n-1}t}{\sin
^{n+1}t} \big[\sin^{n-1}s - (n-1)\Psi _{n-2}(s) \cos s\big]\big[C_n
\Phi _{n-2}(t)\cos t - \sin s\big]
\end{align}
for $s \in (0, t +\pimez)$, where
\begin{equation}\label{defcn}
C_n=\frac{4\omega_{n-1}^2}{n\omega_{n}^2}.
\end{equation}
It is easily seen that
$$ \sin^{n-1}s - (n-1)\Psi _{n-2}(s) \cos s  > 0$$
for $s\in (0, \pi]$. Thus, on setting
\begin{equation}\label{fn}
f_n(t)=C_n\Phi_{n-2}(t)\cos t,
\end{equation}
one has that $\frac{\partial F}{\partial s}(t,s)=0$ if
\begin{equation}
\label{condmv} \sin s=f_n(t).
\end{equation}
We claim that
\begin{equation}
\label{cienne} 0<f_n(t)<1 \qquad \hbox{for $t\in(0,\pimez)$.}
\end{equation}
The first inequality in \eqref{cienne} is trivial. The second
inequality can be established  by induction. If $n=2$, then
$f_2(t)=\frac8{\pi^2}t\cos t$, and an elementary analysis of
$f_2'(t)$ ensures that
\begin{equation*}
\max_{t\in(0,\pimez]}f_2(t)=\max_{t\in(0,\tfrac\pi3]}f_2(t)\le \frac8{3\pi}<1.
\end{equation*}
If $n=3$,  then $f_3(t)=\frac3{4}\sin t\cos t$, whence
\begin{equation*}
f_3(t)\le \frac3{4}<1 \qquad \hbox{for $t\in(0,\pimez)$.}
\end{equation*}
Finally,
\begin{equation*}
f_{n+2}(t)=C_{n+2}\Phi_{n}(t)\cos t\le C_{n+2}\Phi_{n-2}(t)\cos t=\frac{n(n+2)}{(n+1)^2}f_{n}(t)\le f_{n}(t),
\end{equation*}
for every $n\ge2$, where the first inequality holds since $\Phi_{n}
\leq \Phi_{n-2}$. This completes the proof of \eqref{cienne}.
\\
As a consequence of \eqref{cienne},  equation \eqref{condmv} admits
an unique solution $s_t$ in $(0,\pimez]$, given by
\begin{equation}
\label{condmvs} s_t=\arcsin(C_n\Phi_{n-2}(t)\cos t).
\end{equation}
Moreover,  $C_n \Phi _{n-2}(t)\cos t - \sin s <0$, if $s_t<s <
\min\{\pi-s_t,t+\pimez\}$, and  $C_n \Phi _{n-2}(t)\cos t - \sin s
\ge0$ otherwise. Hence,
\begin{align}\label{supF}
\sup_{s\in(0,t+\pimez)} F(t,s)=\max\Big\{F(t,s_t),\lim_{s\rightarrow
(t+\pimez)^-}F(t,s)\Big\}.
\end{align}
Thus, inequality \eqref{ineqmv} will follow if we show that
\begin{align}\label{F1}
\lim_{s\rightarrow (t+\pimez)^-}F(t,s)<0
\end{align}
and $F(t, s_t) < 0$, namely
\begin{align}\label{condmvequiv}
\Phi_{n}(t)\sin s_t+\frac{\cos^n
t}{\sin^{n-1}s_t}(\Psi_{n-2}(s_t)-\Psi_{n}(s_t))-\frac1{n-1}\cos
t>0.
\end{align}
 As far as  inequality \eqref{F1} is concerned,
note that
\begin{align}
\lim_{s\rightarrow (t+\pimez)^-}F(t,s)&=\displaystyle 1-4\frac{n-1}n
\frac{\omega_{n-1}^2}{\omega_{n}^2} \textstyle\Phi_{n-2}(t)\left(\Phi_{n}(t)-\Psi_{n}(t+\pimez)\right)-(n-1)\Psi_{n-2}(t+\pimez)=
\label{F2}\\
\notag\\ \notag &=\Big(1+ \frac{2(n-1)\omega_{n-1}}{n \omega
_n}\Phi_{n-2}(t)\Big)\left(1-\frac{n\omega_n}{2\omega_{n-1}}\right),
\end{align}
where the second equality holds thanks to the fact that
$$\Psi _n(t+\pimez) = \Psi _n(\pimez) + \Phi_n (t) =
\frac{\omega _n}{2 \omega _{n-1}} + \Phi_n (t),$$ and
$$\Psi _{n-2}(t+\pimez) = \Psi _{n-2}(\pimez) + \Phi_{n-2} (t) =
\frac{n \omega _n}{2 (n-1) \omega _{n-1}} + \Phi_{n-2} (t).$$
Observe that
\begin{equation}
\label{aenne} a_n=\frac{n\omega_n}{\omega_{n-1}}> 2 \qquad \hbox{for
$n\ge2$.}
\end{equation}
Inequality \eqref{aenne} follows by induction, from the fact that
 $a_2=\pi>2$, $a_3=4>1$ and
\begin{equation}
\label{A1} a_{n+2}=\frac{n+1}{n} a_n>a_n ,
\end{equation}
for $ n\ge2$, inasmuch as
\begin{align}\label{nov21}
\frac{a_{n+2}}{a_n} = \frac{\frac{(n+2)\omega _{n+2}}{\omega
_{n+1}}}{\frac{n \omega _n}{\omega _{n-1}}} =
\frac{n+2}{n}\frac{\Gamma (1+ \frac n2) \Gamma (1+
\frac{n+1}2)}{\Gamma (1+ \frac{n+2}2) \Gamma (1+ \frac {n-1}2)}=
\frac{n+2}{n} \frac{n+1}{n+2} = \frac{n+1}n.
\end{align}
Equation \eqref{F1} is a consequence of \eqref{F2} and
\eqref{aenne}.
\\ Let us now focus on
 inequality \eqref{condmvequiv}. We begin by showing that
\begin{equation}
\label{condmva}
\Psi_{n-2}(s)-\Psi_{n}(s)\ge\frac1{n-1}\sin^{n-1}s-\frac1{n+2}\sin^{n+1}s\qquad
\hbox{for $s\in[0,\pimez]$.}
\end{equation}
To see this, define $h: [0,\pimez] \to \mathbb R$ as
$$h(s)=\Psi_{n-2}(s)-\Psi_{n}(s)-\frac1{n-1}\sin^{n-1}s+\frac1{n+2}\sin^{n+1}s \quad    \hbox{for $s \in [0,\pimez]$,}$$
and notice that
\begin{equation*}
h'(s)=h_1(s)\sin^{n-2}s\cos s \quad    \hbox{for $s \in
[0,\pimez]$,}
\end{equation*}
where $h_1(s)=\cos s -1+ \frac{n+1}{n+2}\sin^2s$.
 An analysis of the monotonicity properties of  $h_1$ tells us
   that there exists $\bar s\in(0,\pimez)$ such that $h(s)$ is increasing in $[0,\bar s]$ and decreasing in $[\bar s,\pimez]$.
   Therefore, inequality \eqref{condmva} will follow if we show that
\begin{equation}
\label{condmvb}
\Psi_{n-2}(\pimez)-\Psi_{n}(\pimez)-\frac1{n-1}+\frac1{n+2}>0.
\end{equation}
Since
$$\Psi_n(\pimez) = \frac{\omega _n}{2 \omega _{n-1}} \quad \hbox{and} \quad
  \Psi_{n-2}(\pimez)=\frac{n\omega_n}{2(n-1)\omega_{n-1}},
$$
on setting $$b_n=\frac{(n+2)\omega_n}{\omega_{n-1}},$$
inequality
\eqref{condmvb} is equivalent to
\begin{equation}
\label{condmvbn} b_n>6 \qquad \hbox{for $n\ge2$.}
\end{equation}
Inequality \eqref{condmvbn} trivially holds  for $n=2,3$. Also,
$$b_{n+2}=(n+4)\frac{\omega_{n+2}}{\omega_{n+1}}=b_n\frac{(n+4)(n+1)}{(n+2)^2}>b_n
$$
for $n\ge2$. Hence, inequality \eqref{condmvbn} follows by
induction. The  proof of \eqref{condmva} is thereby complete.
\\
On recalling \eqref{condmvs} and \eqref{condmva}, and making use of
\eqref{nov24}, in order to accomplish the prove of inequality
\eqref{condmvequiv}
 it thus suffices to show that
\begin{equation}
\label{condmvsimpl}
\frac{n-1}nC_n\Phi_{n-2}(t)^2-\frac{1-\cos^{n-1}t}{n-1}+C_n\Phi_{n-2}(t)\cos^{n-1}t\left(\frac{\sin
t}n-\frac{C_n}{n+2}\Phi_{n-2}(t)\cos^2t\right)>0.
\end{equation}
Assume first that $n\ge3$, and define the functions
\begin{align}
\quad&k_n(t)=\frac{n-1}nC_n\Phi_{n-2}(t)^2-\frac{1-\cos^{n-1}t}{n-1},%\label{kappa1}
\\
\notag\\
\quad&\kappa_n(t)=\frac{\sin t}n-\frac{C_n}{n+2}\Phi_{n-2}(t)\cos^2t
%\label{kappa2}
\end{align}
for $t\in(0,\pimez)$. Let us first take into account the function
$k_n$.  If $n=3$, then
\begin{equation}\label{febbraio203}
 k_3(t)=0 \qquad \hbox{ for every $t\in(0,\pimez)$.}
 \end{equation}
When $n\ge4$, one has that
\begin{equation*}
%\label{condmvbn}
k_n'(t)=\left(2\frac{n-1}nC_n\Phi_{n-2}(t)-\sin t\right)\cos^{n-2}t.
\end{equation*}
We claim  that there exists $\bar t\in(0,\pimez)$ such that $k_n(t)$
is increasing in $(0,\bar t\,]$ and
 decreasing in $[\bar t,\pimez)$. Our claim follows from
 the fact that $k_n'(0)=0$ and
\begin{equation}\label{febbraio202}
 2\frac{n-1}nC_n\Phi_{n-2}(\pimez)<1.
 \end{equation}
  The latter property is in turn a consequence of the  equality
\begin{equation*}
2\frac{n-1}nC_n\Phi_{n-2}(\pimez)=\frac{4\omega_{n-1}}{n\omega_n}=\frac4{a_n},
\end{equation*}
where $a_n$ is the sequence defined by \eqref{aenne}, owing to
\eqref{A1}, and to  the fact that $a_3=4$ and that $a_4= \frac 83
\pi
>4$, by   the first equality
in \eqref{nov21}. Altogether, since $k_n(0)=k_n(\pimez)=0$, we have
proved that, if $n\ge4$, then
\begin{equation}
\label{kappa1} k_n(t)>0 \qquad \hbox{for $t\in(0,\pimez)$.}
\end{equation}
Consider next the function $\kappa_n$.  One has  that
\begin{equation}\label{nov25}
\kappa_n'(t) = \frac{\cos t}n \left( 1 - \frac{n C_n}{n+2}
\cos^{n-1}t + \frac{2 n C_n}{n+2} \Phi _{n-2} (t) \sin t\right)
\quad \hbox{for $t\in(0,\pimez)$.}
\end{equation}
We claim that
\begin{equation}\label{nov30}
 \frac{n}{n+2}C_n<1 \quad \hbox{for $n \geq 2$.}
 \end{equation}
 To verify inequality \eqref{nov30}, observe that both the subsequence
 $\frac{2n}{2n+2}C_{2n}$ and the subsequence
 $\frac{2n+1}{(2n+1)+2}C_{2n+1}$ are increasing. This is a
 consequence of the inequality $\frac{n}{n+2}C_n <
 \frac{n+2}{n+4}C_{n+2}$, which holds for every $n \geq 2$ and
 follows from the fact that $\Gamma (s+1) = s \Gamma (s)$ for every
$s \in \mathbb R$. On the other hand,  Stirling's formula implies
that
\[\frac{n}{n+2}C_n=\frac{4\omega_{n-1}^2}{(n+2)\omega_{n}^2}\le\lim_n\frac{4\omega_{n-1}^2}{(n+2)\omega_{n}^2}=\pimez<1.
\]
Combining these pieces of information yields inequality
\eqref{nov30}. Owing to this inequality, we infer from \eqref{nov25}
that $\kappa_n'(t)>0$ for $t\in(0,\pimez)$. Hence, since
$\kappa_n(0)=0$,
\begin{equation}
\label{kappa2} \kappa_n(t)>0\qquad \hbox{for $t\in(0,\pimez)$.}
\end{equation}
Coupling either \eqref{febbraio203} or \eqref{kappa1}, with
\eqref{kappa2} yields \eqref{condmvsimpl}, and hence inequality
\eqref{condmvequiv} for $n\ge3$.
\par
Let us finally consider the case when $n=2$. Observe that, in this
case, the left-hand side of \eqref{condmvsimpl} equals
$$\tilde k + C_2 \Phi _0 (t) \tilde \kappa (t)\cos t \qquad \hbox{for $t\in(0,\pimez)$,}$$
 where we have set
\begin{align}
\quad&\tilde k(t)=k_2(t)+\frac{16}{\pi^4}t^2\cos t=\frac4{\pi^4}t^2(\pi^2+4\cos t)-1+\cos t,\label{tkappa1}\\
\notag\\
\quad&\tilde\kappa(t)=\kappa_2(t)-\frac
2{\pi ^2} t=\frac1{2}\sin t
-\frac{2}{\pi^2}t(1+\cos^2t).\label{tkappa2}
\end{align}
We claim that
\begin{align}
\quad& m_1(t)\equiv t-\frac{\pi^2}2\sqrt{\frac{1-\cos t}{\pi^2+4\cos t}}>0& \hbox{for $t\in(0,\pimez)$},\label{ausil1}\\
\notag\\
\quad&m_2(t)\equiv \frac{\sin t}{1+\cos^2t}-\frac4{\pi^2}t>0& \qquad
\hbox{for $t\in(0,\pimez)$}.\label{ausil2}
\end{align}
Inequalities \eqref{ausil1} and \eqref{ausil2} imply $\tilde k(t)>0$
and $\tilde \kappa(t)>0$, and hence  \eqref{condmvsimpl} follows,
thus establishing inequality \eqref{condmvequiv}  also for $n=2$.
%and a close inspection of the above arguments also gives the equality case in \eqref{ineqmv}.
It just remains to prove inequalities \eqref{ausil1} and
\eqref{ausil2}. An analysis  of monotonicity properties of the
function $m_1(t)$ tells us that there exists $\tilde t\in (0,
\pimez)$ such that $m_1(t)$ is increasing in $(0,\tilde t\,]$ and
decreasing in $[\tilde t, \pimez)$. Inequality \eqref{ausil1} hence
follows, since $m_1(0^+)=m_1( \pimez ^-)=0$. \\ As for
 \eqref{ausil2},   a study of the sign of
 $m_2''(t)$ tells us that
 there exists $\hat t\in (0, \pimez)$ such that $m_2(t)$ is convex in $(0,\hat t\,]$, and concave in $[\hat t,\pimez)$.
  This piece of information, combined with the fact that $m_2(0^+)=0$, $m_2'(0^+)=\frac{\pi^2-8}{2\pi^2}>0$ and $m_2( \pimez ^-)=\frac{\pi-2}{\pi}>0$,
  yields \eqref{ausil2}.
\qed

\bigskip
\par
%
%Before proving Theorem \ref{isopconstmv} we recall the definition of
%spherical symmetrization and an approximation result via polyhedra.

Given a measurable set $E \subset
\B$,  we denote by $E^\sharp$ the spherical symmetral of $E$ about
the half-axis $X=\{(x_1, \dots , x_n): x_1\geq 0, x_2=\dots = x_n
=0\}$.  The set $E^\sharp$ is defined as the subset of $\B$ such that
the intersection of $E^\sharp$ with any sphere $S$ centered at $0$
is a spherical cap, centered at $S \cap X$, such that $\hh (E^\sharp
\cap S) = \hh (E \cap S)$. In particular, $E^\sharp$ is symmetric
about the $x_1$-axis.
\\ The very definition of spherical symmetrization, and the use of polar coordinates, ensure that
\begin{equation}\label{febbraio204}
|E| = |E^\sharp|
\end{equation}
for every measurable subset $E$ of $\B$. The definition of spherical
symmetrization again tells us that
\begin{equation}\label{iso4}
\hh ( \partial E^\sharp \cap \partial \B ) = \hh (\partial E \cap
\partial \B)
\end{equation}
if $E$ is a sufficiently regular subset of $\B$. Moreover, a
classical property of spherical symmetrization entails that it does
not increase
 perimeter relative to $\B$ of regular subsets $E$ of $\B$; namely
\begin{equation}\label{iso5}
\hh (\partial  E^\sharp \cap  \B) \leq  \hh (\partial  E \cap
 \B),
\end{equation}
see e.g. \cite{Ka1}. In fact, equations \eqref{iso4} and
\eqref{iso5} will be exploited when $E$ is just the intersection of
$\B$ with a polyhedron. This will suffice for our purposes, since we
shall make use of a result from \cite[Lemma 9.4.1/3]{Mazbook} which
tells us that, given any measurable set  $E \subset \B $   such that
$ \hh (\partial ^M E ) < \infty$, there exists a sequence of
polyhedra $\{P_k\}$ in $\rn$ with the following properties. Define
$Q_k= P_k \cap \B$ for $k \in \N$. Then
\begin{equation}\label{iso2'}
\lim _{k \to \infty} \chi _{Q_k} = \chi _E
\end{equation}
 in
$L^1(\B)$,
\begin{equation}\label{iso2} \lim _{k \to
\infty}  \hh (\partial  Q_k \cap \B) = \hh (\partial ^M E \cap \B),
\end{equation}
and
\begin{equation}\label{iso3} \lim _{k \to
\infty}  \hh (\partial Q_k \cap \partial \B) = \hh (\partial ^M E
\cap
\partial \B).
\end{equation}

\medskip
\par
\noindent{\bf Proof of Theorem \ref{isopconstmv}}. Let ${\cal
Q}_{\rm mv}$ be the functional defined by \eqref{defQmv}. Since for
any measurable set $E \subset \B$ such that $ \hh (\partial ^M E ) <
\infty$ there exists a sequence of polyhedra $\{P_k\}$ satisfying
\eqref{iso2'}--\eqref{iso3}, one has that
\begin{align}\label{nov35}
\sup _{E\subset \B}{\cal Q}_{\rm mv}(E) & = \sup _{
\begin{tiny}
 \begin{array}{c}{Q=P \cap \B} \\
P \hbox{ is a polyhedron}
 \end{array}
  \end{tiny}
} {\cal Q}_{\rm mv}(Q)
\leq \sup _{
\begin{tiny}
 \begin{array}{c}{Q=P \cap \B} \\
P \hbox{ is a polyhedron}
 \end{array}
  \end{tiny}
} {\cal Q}_{\rm mv}(Q^\sharp)
\\ \nonumber & \leq  \sup
 _{
\begin{tiny}
 \begin{array}{c}{E\subset \B} \\
E= E^\sharp
 \end{array}
  \end{tiny}
} {\cal Q}_{\rm mv}(E)
  \leq
\sup _{E\subset \B}{\cal Q}_{\rm mv}(E).
\end{align}
Note that the first equality in \eqref{nov35} holds by
\eqref{iso2'}--\eqref{iso3}, and the first inequality by
\eqref{febbraio204} and \eqref{iso4}.
 Hence,
$$
K_{\rm mv}(\B)= \sup _{E\subset \B}{\cal Q}_{\rm mv}(E) = \sup
 _{
\begin{tiny}
 \begin{array}{c}{E\subset \B} \\
E= E^\sharp
 \end{array}
  \end{tiny}
} {\cal Q}_{\rm mv}(E).
$$
Thus, we may limit ourselves to maximize ${\cal Q}_{\rm mv}(E)$ in
the class of sets  $E$ such that $E= E^\sharp$, and hence, in
particular, $\partial ^M E \cap\partial \B$ is a spherical cap (or
an empty set) on $\partial \B$. Since the functional ${\cal Q}_{\rm
mv}(E)$ is invariant under replacements of $E$ with $\B \setminus
E$, we may also assume that $\hh(\partial ^M E \cap\partial \B)\ge
\hh(\partial \B\setminus\partial ^M E)$. Let us also observe that
${\cal Q}_{\rm mv}(E)$ cannot achieve its maximum at any set  $E$
such that $\hh(\partial \B\setminus\partial ^M E)=0$. Indeed, if
this equality holds, then $\hh (\partial ^M E \cap \B) = \hh
(\partial ^M (\B\setminus E) \cap \B) = \hh (\partial ^M
(\B\setminus E))$, and hence
\begin{equation}\label{*}
{\cal Q}_{\rm mv}(E)=\frac{|\B\setminus E|\hh (
\partial \B)}{|\B|\hh (\partial ^M (\B\setminus
E))}\le1<\frac{n\omega_n}{2\omega_{n-1}}.
\end{equation}
Observe that the last inequality in \eqref{*} holds by
\eqref{aenne}. The first inequality is instead a consequence of the
fact that, by the standard isoperimetric theorem, the ratio
$\frac{|\B\setminus E|}{\hh (\partial ^M (\B\setminus E))}$ does not
decrease if $\B\setminus E$ is replaced with a ball of equal
Lebesgue measure, and that it increases if the ball is replaced with
a larger ball. On the other hand, the rightmost side of \eqref{*}
agrees with the functional ${\cal Q}_{\rm mv}$ evaluated at a
half-ball. Altogether,
\begin{equation}\label{march200}
\sup _{E\subset \B}{\cal Q}_{\rm mv}(E) = \sup _{E\in \mathcal
E}{\cal Q}_{\rm mv}(E)\,,
\end{equation}
where
\begin{equation}\label{march201}
\mathcal E =\{E \subset \B : E \, \hbox{is measurable}\, , E =
E^\sharp , \, \hh(\partial ^M E \cap\partial \B)\ge \hh(\partial
\B\setminus\partial ^M E)>0\}.
\end{equation}
\iffalse

we have used isoperimetric inequality and \eqref{aenne}, and we
recall that the last term corresponds to ${\cal Q}_{\rm mv}$
evaluated on a half-ball. So we assume
\begin{equation}\label{condacca}
\hh(\partial \B\setminus\partial ^M E)>0.
\end{equation}

Let us denote by $\mathcal E$ the class of subsets of $\B$
fulfilling these properties. Fix any set $E \in \mathcal E$. One has
that
\begin{equation}\label{Qmvalt}
{\cal Q}_{\rm mv}(E)=\frac{|E| [\hh (\partial \B \setminus \partial
^M E)-\hh (\partial ^M E \cap
\partial \B)] + |\B|\hh (\partial ^M E \cap
\partial \B)}{|\B|\hh (\partial ^M E \cap \B)},
\end{equation}
where the expression in square brackets is non-positive.

\fi
 Given $t \in \mathbb R$, define the  half-space
\begin{equation}\label{Ht}
H_t=\{x\in\R^n:x_1>t\},
\end{equation}
%{\color{blue}For any $E \in
%\mathcal E$ satisfying \eqref{condacca},}
and set $\bar t=\sup\{t\in\R:H_t\supset\partial ^M E \cap
\partial \B\}$. Then either $|H_{\bar t}\cap \B|>|E|$, or $|H_{\bar t}\cap \B|\le|E|$.
\\
Assume first that $|H_{\bar t}\cap \B|>|E|$, and
  consider a ball $ B$ such that
$ B\cap \partial \B=\partial H_{\bar t} \cap
\partial \B$ and $|\B\setminus  B|=|E|$.
Define
$$\widetilde E =
(\B\backslash E)\cup( B\backslash  \B).$$
 Clearly,
$$|\widetilde E| = | B|,$$
and, by the isoperimetric property of the ball,
$$\hh (\partial \widetilde E) \geq \hh (\partial  B).$$
On the other hand,
$$\hh (\partial \widetilde E) = \hh (\partial  B \cap (\rn \setminus \B)) +
\hh (\partial ^M E \cap \B),$$ and
$$\hh (\partial  B) = \hh (\partial  B \cap (\rn \setminus \B)) +
\hh (\partial  B \cap \B).$$ Hence,
$$\hh(\partial B\cap \B)\le\hh(\partial ^M E \cap
 \B).
$$
On setting  $\widehat E=\B\backslash  B$, one has that
 $|\widehat E|=|E|$ and
$${\cal Q}_{\rm mv}(E)\le{\cal Q}_{\rm mv}(\widehat E).
$$
Suppose next that $|H_{\bar t}\cap \B|\le|E|$. Then  the set
$\widehat E=H_{\bar t}\cap \B$ satisfies the inequalities $|\widehat
E|\le|E|$, and $\hh (\partial E \cap \B) \geq \hh (\partial H_{\bar
t} \cap \B) = \hh (\partial \widehat E \cap \B)$, whence
$${\cal Q}_{\rm mv}(E)\le{\cal Q}_{\rm mv}(\widehat E).
$$
Altogether, we have shown that
\begin{equation}\label{Klemma}
K_{\rm mv}(\B)= \sup_{ \scriptstyle E\in{\cal A}}{\cal Q}_{\rm
mv}(E),
\end{equation}
where ${\cal A}$ denotes the collection of those subsets $E$ of $\B$
such that  $E$  is the complement in $\B$ of either a ball, or of a
half-space, with $ \hh(\partial ^M E \cap\partial \B)\ge
\hh(\partial \B\setminus\partial ^M E)>0$.
%\begin{align}
%& E \text{ is the complement of a ball or of a half-space in }\B;\label{AA}\\ \notag\\
%& \hh(\partial ^M E \cap\partial \B)\ge
%\hh(\partial \B\setminus\partial ^M E).\label{AA1}
%\end{align}
\\ In view of \eqref{Klemma}, in order to conclude our proof it
remains to show that
\begin{equation}\label{mvequiv}
\sup_{ \scriptstyle E\in{\cal A}}{\cal Q}_{\rm mv}(E)\le
\frac{n\omega_n}{2\omega_{n-1}}.
\end{equation}
%To this purpose, observe first that, if $E=\B\setminus B_r$ for some
%ball $B_r\subset \B$   of radius $r$, $0<r<1$, then ${\cal Q}_{\rm
%mv}(E)=r<1$, and hence strict inequality holds in \eqref{mvequiv},
%owing to \eqref{aenne}.
%\begin{equation}\label{mveasy}
%{\cal Q}_{\rm mv}(E)=r<1
%\end{equation}
\\ We may thus focus on
 the case when $E=E_{\vartheta ,
\varphi}$ for some $(\vartheta,\varphi)\in\Theta$, where $\Theta$ is
defined in \eqref{defdelta}. In other words, we have to show that
\begin{equation}\label{mvcompl}
{\cal Q}_{\rm mv}(E_{\vartheta ,
\varphi})\le\frac{n\omega_n}{2\omega_{n-1}},\qquad(\vartheta,\varphi)\in\Theta,
\end{equation}
 the equality being attained  if
$(\vartheta , \varphi) = (\pimez , 0)$, in which case   $E_{\pimez ,
0}$ is a half-ball.
\\
%Since the quantity ${\cal Q}_{\rm mv}(E_{\vartheta , \varphi})$ is
%given by
Owing to formula \eqref{QEmv}, the conclusion follows from Lemma
\ref{symmmv1}.
 \qed

\section{A trace inequality on $\mathbb B^n$ with median normalization}\label{sec1}

Our concern in this section is to detect the extremal functions in
the Poincar\'e trace inequality
\begin{equation}\label{Cmedoptn}
\|\widetilde u - {\rm med}_{\B,\sigma} (u)\|_{L^1(\partial \B )}
\leq C_{\rm med}(\B, \sigma)\|Du\|(\B),
\end{equation}
with optimal constant $C_{\rm med}(\B, \sigma)$, for $u \in BV(\B)$.
Recall that ${\rm med}_{\B,\sigma} (u)$ denotes the $\sigma$-median
of $u$, defined as in \eqref{medsigma}, which agrees with the usual
median ${\rm med}_{\B} (u)$ when $\sigma = \tfrac 12$.
\\ This is the content of the following theorem.

\begin{theorem}\label{isopconstpmed} Let $n \geq 2$ and let $\sigma\in(0,1)$. Set $\rho = \max\{\sigma, 1-\sigma\}$.
Then equality holds in \eqref{Cmedoptn} if $u$ is the characteristic
function of the half-moon shaped set
 $E_{\vartheta_\rho,\varphi_\rho}$   as in Figure
\ref{intersection},   where $(\vartheta_\rho, \varphi_\rho)$ is the
unique solution in the set $\Upsilon$ (defined by \eqref{dominio})
to the system
\begin{equation}\label{sistema}
\quad\left\{
\begin{array}{l}
\displaystyle\frac{\Psi_{n-2}(\varphi)}{\Psi_{n-2}(\vartheta)}\frac{\sin^{n}\vartheta}{\sin^{n}\varphi}
=1-{\frac{\rho(n-1)\Psi_{n-2}(\pi)\cos\vartheta}{(n-1)\cos
\vartheta\,\Psi_{n-2}({\vartheta})-\sin^{n-1}\vartheta}}
\\
\\
\displaystyle\frac{\cos\varphi}{\sin\varphi}=\frac{\cos\vartheta}{\sin\vartheta}\left(1-{\frac{\rho(n-1)\Psi_{n-2}(\pi)}{(n-1)\cos^2
\vartheta\,\Psi_{n-2}({\vartheta})-\sin^{n-1}\vartheta\cos\vartheta}}\right).
\end{array}
\right.
\end{equation}
\end{theorem}

\medskip
\par

 Theorem \ref{isopconstpmed} follows from Theorem \ref{Comegasigma},
via the  next result.

\begin{figure}\setlength\fboxsep{30pt} \setlength\fboxrule{0.0pt}
\hskip3.7cm\fbox{\includegraphics[width=2.2 in]{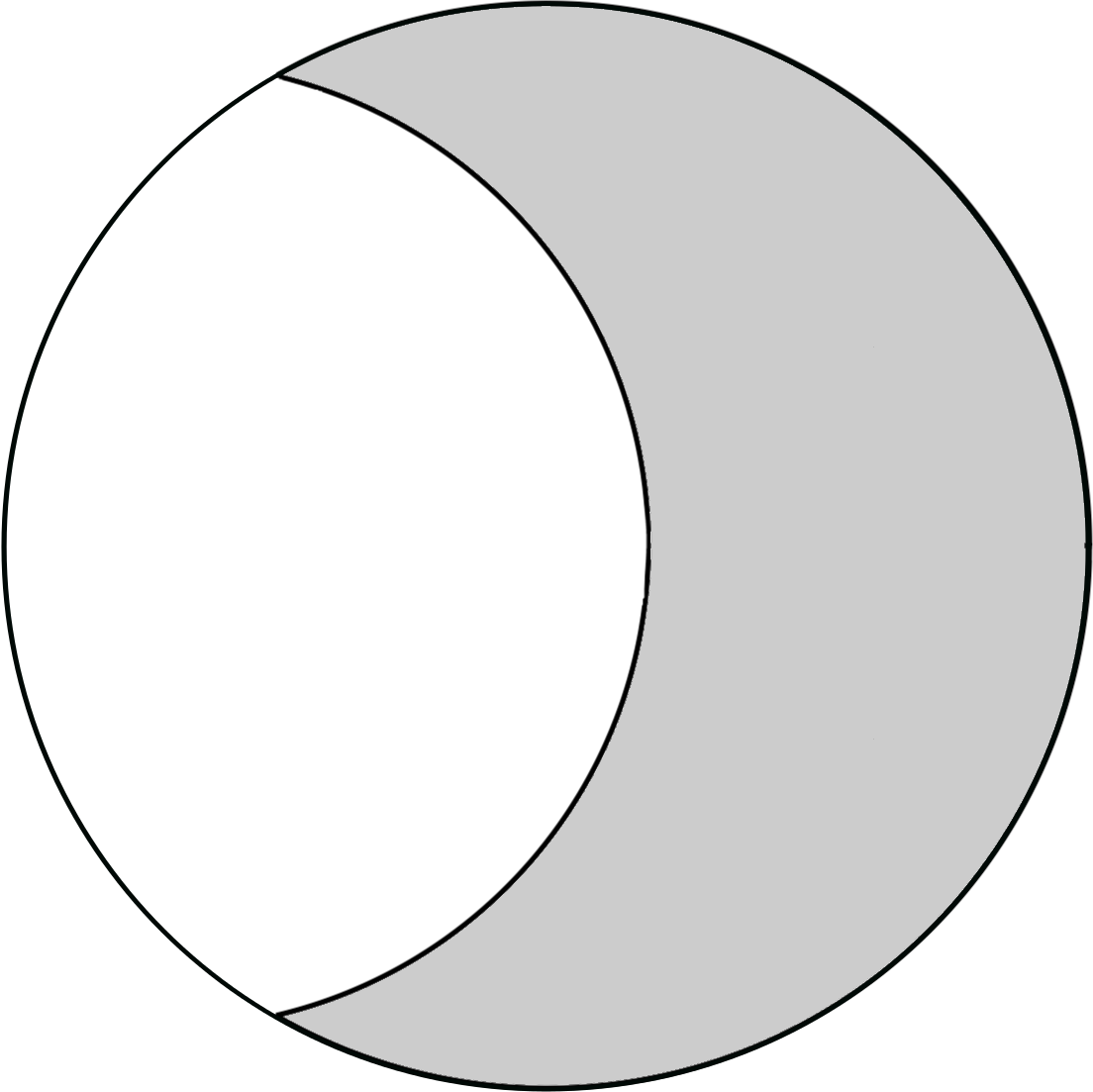}}
\caption{A half-moon shaped extremal set for $K_{\rm med}(\B,
\sigma)$}\label{intersection_bis}
\begin{picture}(0,0)(150,-30)
\put(415,180){$\B$}
\put(410,95){$E$}
\end{picture}
\end{figure}
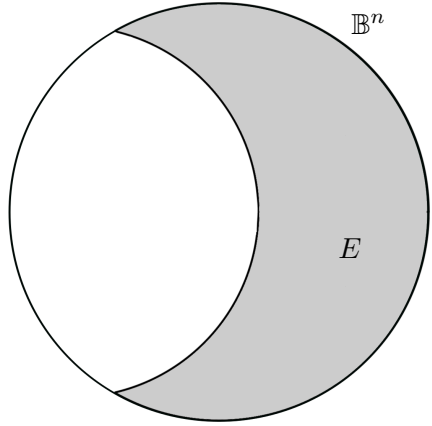

\begin{theorem}\label{isopconst} Let $n \geq 2$ and let $\sigma\in(0,1)$. Then the set
$E_{\vartheta_\sigma,\varphi_\sigma}$, defined as in Theorem
\ref{isopconstpmed} with $\rho$ replaced with $\sigma$, is extremal
for $K_{\rm med}(\B, \sigma)$ (Figure \ref{intersection_bis}).
\end{theorem}

Theorem \ref{isopconst} is in turn a straightforward consequence of
Lemmas \ref{symm} and \ref{lune} below. The former enables us to
reduce the detection of extremals for $K_{\rm med}(\B, \sigma)$ to
the class of sets of the form $E_{\vartheta,\varphi}$ with
$(\vartheta_\rho, \varphi_\rho) \in \Upsilon$. The latter identifies
$E_{\vartheta_\sigma,\varphi_\sigma}$, with $(\vartheta_\sigma,
\varphi_\sigma)$  solving \eqref{sistema} with $\rho$ replaced with
$\sigma$, as the unique extremal for $K_{\rm med}(\B, \sigma)$ in
this special class.

In what follows, we  denote  by ${\cal Q}_{\rm med}(E)$ the
functional of $E$ that is maximized on the right-hand side of
\eqref{Komegasigma}, namely
\begin{equation}\label{QEgen} {\cal Q}_{\rm med}(E) =\frac{ \hh (\partial ^M E \cap
\partial \B) }{\hh (\partial ^M E \cap \B)}.
\end{equation}
We also set
\begin{equation}\label{B}
{\cal B} = \{E \subset \B: E = \B \setminus B, \hbox{$B$ is a ball
or a half-space}\,   \hh(\partial \B\setminus\partial ^M E)>0\}.
\end{equation}

%
%A preliminary analysis will enable us to restrict the detection of
%extremal sets for \eqref{Komegasigma} to the class of subsets of
%$\B$ of the form $E_{\vartheta,\varphi}$, with $(\vartheta,\varphi)$
%from the set
%\begin{equation}\label{defdelta_bis}
%\Xi=\{(\vartheta,\varphi)\in \R^2:0< \vartheta<\pi ,\>0\le \varphi<
%\vartheta\}.
%\end{equation}
% Note that, by \eqref{def1}, \eqref{def2},
%\begin{equation}\label{QE} {\cal Q}_{\rm med}(E_{\vartheta,\varphi}) =
%\frac{
%\Psi_{n-2}(\vartheta)\sin^{n-1}\varphi}{\Psi_{n-2}(\varphi)\sin^{n-1}\vartheta}
%\quad \hbox{for $(\vartheta,\varphi) \in \Xi$,}
%\end{equation}
%where, as usual, the right-hand side is extended by continuity for
%$\varphi=0$.

\iffalse
 Observe that the limit case $\varphi=0$ corresponds to the
case where $E_{\vartheta, 0}$ is the complement of a half-space in
$\B$. The ratio \eqref{QEgen} has to be considered also in this case
and we extend by continuity the definition \eqref{defG} for
$\varphi= 0$. Actually, as in the previous section, also the case
where $E$ is the complement of a ball contained in $\B$ has to be
considered. It corresponds to $\vartheta=\pi$, but then only the
couple $(\vartheta,\varphi)=(\pi,\pi)$ is possible. A direct
inspection of the ratio \eqref{QEgen} proves that this case does not
give a isoperimetric set.

\fi

\begin{lemma}\label{symm}
Let $ n \geq 2$, and let $\sigma \in (0, 1)$. Then,
\begin{equation}\label{iso0}\sup _{|E|\leq \sigma |\B|} {\cal Q}_{\rm med}(E)
= \sup_{
\begin{tiny}
 \begin{array}{c}{E\in\cal{B}} \\
|E|= \sigma |\B|
 \end{array}
  \end{tiny}
} {\cal Q}_{\rm med}(E),
\end{equation}
where ${\cal B}$ is defined by \eqref{B}.
%denotes the collection of those subsets $E$ of
%$\,\B$ such that  $E$  is the complement in $\B$ of a ball or of a
%half-space, with $ \hh(\partial \B\setminus\partial ^M E)>0$.
\end{lemma}
{\bf Proof}. Given any measurable set $E \subset \B$ of finite
perimeter in $\B$,  consider a sequence of polyhedra $\{P_k\}$
satisfying \eqref{iso2'}--\eqref{iso3}. Fix any $k \in \N$. By
properties \eqref{iso4} and \eqref{iso5}  of spherical
symmetrization,
 for every $\varepsilon >0$ and $\delta \in (0, 1-\sigma )$,
there exists $\overline k \in \N$ such that, if $k \geq \overline
k$, then
\begin{equation}\label{iso6}
|Q_k^\sharp| \leq (1+\delta )|E|\,,
\end{equation}
and
\begin{equation}\label{iso7}
 {\cal Q}_{\rm med}(E)
 \leq \frac{ \hh (\partial Q_k^\sharp\cap
\partial \B)+\varepsilon}{\hh (\partial Q_k^\sharp \cap \B)-\varepsilon}.
\end{equation}
Owing to the arbitrariness of $\varepsilon$, inequalities
\eqref{iso6} and \eqref{iso7} imply that
\begin{equation}\label{iso8}
\sup _{|E|\leq \sigma |\B|}  {\cal Q}_{\rm med}(E) \leq \sup_{
\begin{tiny}
 \begin{array}{c}{E=E^\sharp}\\
|E|\leq (\sigma + \delta) |\B|
 \end{array}
  \end{tiny}
}  {\cal Q}_{\rm med}(E).
\end{equation}
We next show that, for every $\delta$ as above,
\begin{equation}\label{iso9}
\sup_{
\begin{tiny}
 \begin{array}{c}{E=E^\sharp}\\
|E|\leq (\sigma + \delta) |\B|
 \end{array}
  \end{tiny}
}  {\cal Q}_{\rm med}(E)
 =
\sup_{
\begin{tiny}
 \begin{array}{c}{ E\in\cal B} \\
|E|\leq (\sigma + \delta) |\B|
 \end{array}
  \end{tiny}
}
 {\cal Q}_{\rm med}(E).
\end{equation}
An  argument analogous to that in the proof of \eqref{*} tells us
that, if $E$ is any subset of $\B$ such that $\hh(\partial
\B\setminus\partial ^M E)=0$, then
 %Let us observe that any set $E$ such that $\hh(\partial
%\B\setminus\partial ^M E)=0$ cannot achieve the maximum of ${\cal
%Q}_{\rm med}(E)$. Indeed, in such a case, we have:
\begin{equation*}
{\cal Q}_{\rm med}(E)=\frac{ \hh (\partial \B)}{\hh (\partial
(\B\setminus E))}\le1<\frac{n\omega_n}{2\omega_{n-1}}.
\end{equation*}
Moreover, the last expression agree with the functional ${\cal
Q}_{\rm med}$ evaluated at a half-ball in $\B$. We may thus assume
that the sets $E$ on the right-hand side of \eqref{iso9} fulfil the
condition
%
%
%where we have used isoperimetric inequality and \eqref{aenne}, and
%we recall that the last term corresponds to ${\cal Q}_{\rm med}$ evaluated on a half-ball. So we assume
\begin{equation}\label{condaccab}
\hh(\partial
\B\setminus\partial ^M E)>0.
\end{equation}
Let $E$ be subset of $\B$ such that $E=E^\sharp$, $|E|\leq (\sigma +
\delta)|\B|$ and \eqref{condaccab} holds. Define $\bar
t=\sup\{t\in\R:H_t\supset\partial ^M E \cap
\partial \B\}$, where $H_t$ stands for the half-space introduced in the proof of
Theorem \ref{isopconstmv}. Assume first that $|H_{\bar t}\cap
\B|>|E|$. Consider a ball $\widehat B$ such that $\widehat B\cap
\partial \B=\partial H_{\bar t} \cap
\partial \B$ and $|\B\backslash \widehat B|=|E|$.
Set
$$\widetilde E =
(\B\backslash E)\cup(\widehat B\backslash  \B).$$ Clearly,
$$|\widetilde E| = |\widehat B|,$$
and, by the isoperimetric property of the ball,
$$\hh (\partial \widetilde E) \geq \hh (\partial \widehat B).$$
On the other hand,
$$\hh (\partial \widetilde E) = \hh (\partial \widehat B \cap (\rn \setminus \B)) +
\hh (\partial ^M E \cap \B),$$ and
$$\hh (\partial \widehat B) = \hh (\partial \widehat B \cap (\rn \setminus \B)) +
\hh (\partial \widehat B \cap \B).$$
Hence,
$$\hh(\partial\widehat B\cap \B)\le\hh(\partial ^M E \cap
 \B).
$$
Now, if we define  $\widehat E=\B\backslash \widehat B$, then
 $|\widehat E|=|E|$ and
\begin{equation}\label{febbraio300}
{\cal Q}_{\rm med}(E)\le{\cal Q}_{\rm med}(\widehat E).
\end{equation}
In the case when $|H_{\bar t}\cap B|\le|E|$,  the set $\widehat E$
defined as $\widehat E=H_{\bar t}\cap \B$ has the property that
$|\widehat E|\le|E|$ and inequality \eqref{febbraio300} still holds.
 Altogether, equation \eqref{iso9} is  established. From \eqref{iso8}
and \eqref{iso9} we deduce that
\begin{equation}\label{nov1}\sup _{|E|\leq \sigma |\B|} {\cal Q}_{\rm med}(E)
\leq \sup_{
\begin{tiny}
 \begin{array}{c}{ E\in\cal B}\\
|E|\leq (\sigma + \delta) |\B|
 \end{array}
  \end{tiny}
} {\cal Q}_{\rm med}(E).
\end{equation}
In order to accomplish the proof of \eqref{iso0}, it remains to show
that
 \eqref{nov1} continues to hold with $\delta =0$. To verify this
 fact, choose $\delta = \frac 1k$, with $k \in \mathbb N$, in
 \eqref{nov1}, and denote by $
E_k$ a set from $\cal B$ such that $|E_k| \leq \sigma + \frac 1k$
and
\begin{equation}\label{nov2}
{\cal Q}_{\rm med}(E_k)
 \geq
\sup_{
\begin{tiny}
 \begin{array}{c}{E\in\cal B}\\
|E|\leq (\sigma + 1/k) |\B|
 \end{array}
  \end{tiny}
} {\cal Q}_{\rm med}(E) - \frac 1k.
\end{equation}
Thus,
\begin{equation}\label{nov3}
\sup _{|E|\leq \sigma |\B|} {\cal Q}_{\rm med}(E) \leq {\cal Q}_{\rm
med}(E_k) + \frac 1k.
\end{equation}
If there exist infinitely many values of $k$ such that $|E_k| \leq
\sigma $, then \eqref{iso0} immediately follows from \eqref{nov3}.
If, on the contrary, $\sigma < |E_k|\leq \sigma + \frac 1k $ for
all, but finitely many values of $k$, then there exists a
subsequence $ E_{k_j}$ and a set $ E\in\cal B$ such that $ E_{k_j}
\to  E$, and $|E|= \sigma$, and  \eqref{iso0} follows also in this
case. \qed

Now, observe that, by \eqref{def1} and \eqref{def2},
\begin{equation}\label{QE} {\cal Q}_{\rm med}(E_{\vartheta,\varphi}) =
\frac{
\Psi_{n-2}(\vartheta)\sin^{n-1}\varphi}{\Psi_{n-2}(\varphi)\sin^{n-1}\vartheta}
\quad \hbox{for $(\vartheta,\varphi) \in \Upsilon$,}
\end{equation}
where, as usual, the function on right-hand side is extended by
continuity for $\varphi=0$. Let us denote by $G: \Upsilon \to [0,
\infty)$ this function, namely
\begin{equation}\label{G}
G(\vartheta,\varphi) = \frac{
\Psi_{n-2}(\vartheta)\sin^{n-1}\varphi}{\Psi_{n-2}(\varphi)\sin^{n-1}\vartheta}
\quad \hbox{for $(\vartheta,\varphi) \in \Upsilon$.}
\end{equation}
Also, define, for $\sigma \in (0, 1)$,
\begin{equation}\label{xi}
\Xi (\sigma) = \{(\vartheta,\varphi) \in \Upsilon : \,
|E_{\vartheta,\varphi}| \leq \sigma |\B|\}.
\end{equation}
%and
%\begin{equation}\label{M}
%M(\sigma ) = \sup _{(\vartheta,\varphi) \in \Xi (\sigma)}
%G(\vartheta,\varphi).
%\end{equation}

\begin{lemma}\label{lune} Let $n \geq 2$, and let $\sigma \in (0, 1)$. Then   system
\eqref{sistema},  with $\rho$ replaced with $\sigma$, has a unique
solution $(\vartheta_\sigma,\varphi_\sigma) \in \Upsilon$, and
\begin{equation}\label{march202}
\max _{(\vartheta,\varphi) \in \Xi}{\cal Q}_{\rm
med}(E_{\vartheta,\varphi}) = {\cal Q}_{\rm
med}(E_{\vartheta_\sigma,\varphi_\sigma}).
\end{equation}
\end{lemma}
{\bf Proof}. Let $\beta \in (0, \sigma]$. Equation \eqref{def3}
entails that
\begin{equation}\label{condizione}
|E_{\vartheta,\varphi}|=\beta|\B|
\end{equation}
if and only if
\begin{equation}\label{cond}
\frac{ \Psi_{n}(\vartheta)}{\sin^{n}\vartheta}-\beta \frac{
\Psi_{n}(\pi)}{\sin^{n}\vartheta}-\frac{
\Psi_{n}(\varphi)}{\sin^{n}\varphi}=0.
\end{equation}
%The case when $B_1$ is a half-space, that is, when $E$ is a
%spherical segment, corresponds to the limit case $\varphi=0$. In
%such a case the quotient ${\cal Q}_{\rm med}(E)$ is given by
In the borderline case when $\varphi=0$, which corresponds to a set
$E_{\vartheta , 0}$ obtained as the intersection of $\B$ with a
half-space, one has that
\begin{equation}\label{hsp}
G(\vartheta , 0)=(n-1)\frac{ \Psi_{n-2}(\vartheta
(\beta))}{\sin^{n-1}\vartheta (\beta)},
\end{equation}
where $\vartheta (\beta)\in (0,\pi)$ obeys
\begin{equation}\label{thetabar}
\Psi_{n}(\vartheta (\beta))=\beta \Psi_{n}(\pi).
\end{equation}
Observe that the function on the right-hand side of  \eqref{hsp} is
strictly increasing with respect to $\vartheta(\beta)$, and the
latter is a strictly increasing function of $\beta$. Hence the
maximum  of {$G(\vartheta , 0)$}
under the constraint \eqref{condizione}, with $\beta \in (0,
\sigma]$, is achieved for $\beta = \sigma$.
\\ Similarly, the function $G$ cannot attain its maximum at a couple $(\vartheta , \varphi)$ with $\varphi >0$, unless   condition
 \eqref{condizione} is fulfilled with $\beta = \sigma$.  This is
 verified on recalling the geometric meaning of the function $G$.
 Indeed, assume, by contradiction, that $G$ attains its maximum at
 some point $(\vartheta , \varphi)$ with $\varphi >0$.
Then the corresponding set $E_{\vartheta , \varphi}$ is the
complement
 in $\B$ of some ball. Let $H_t$ be the half-space defined as in \eqref{Ht}
 for $\in \mathbb R$,  set $E_t=(E\cup H_t)\cap \B$, and
$$\widehat t = \inf \{t: H_t \cap B
 \subset E_{\vartheta , \varphi}\}.$$
Then there exists    $t<\widehat t$ such that still
$|E_t|<\sigma|\B|$, but ${\cal Q}_{\rm med}(E_t)>{\cal Q}_{\rm
med}(E_{\vartheta , \varphi})$,
inasmuch as $\hh (\partial E_t \cap
\B) < \hh (\partial E_{\vartheta , \varphi} \cap \B)$.
\\
In view of the above consideration, the maximum on the left-hand
side of \eqref{march202} agrees with the maximum of the function $G$
on the set $\Upsilon$ under the constraint
\begin{equation}\label{conda}
\frac{ \Psi_{n}(\vartheta)}{\sin^{n}\vartheta}-\sigma \frac{
\Psi_{n}(\pi)}{\sin^{n}\vartheta}-\frac{
\Psi_{n}(\varphi)}{\sin^{n}\varphi}=0.
\end{equation}
Note that, if $(\vartheta , \varphi)$  fulfils equation
\eqref{conda}, then $\vartheta \geq \vartheta (\sigma )$. Actually,
if $\varphi =0$, then $\vartheta =\vartheta (\sigma )$. On the other
hand, in the case when $\varphi
>0$, one has that $|E_{\vartheta (\sigma ), \varphi}| <
|E_{\vartheta (\sigma ), 0}| = \sigma |\B|$. Consequently,
$\vartheta > \vartheta (\sigma )$, since \eqref{conda} is equivalent
to $|E_{\vartheta (\sigma ), \varphi}|  = \sigma |\B|$. Equation
\eqref{conda} implicitly defines  $\vartheta$ as  a function of
$\varphi$. Indeed, using the notation $\vartheta(\beta)$ introduced
in \eqref{thetabar}, both the function $\xi_\sigma :
[\vartheta(\sigma), \pi) \to [0,\infty)$, given by
\begin{equation}
\label{TH} \xi_\sigma (\vartheta)=\frac{
\Psi_{n}(\vartheta)}{\sin^{n}\vartheta}-\sigma \frac{
\Psi_{n}(\pi)}{\sin^{n}\vartheta}\quad \hbox{for $\vartheta \in
[\vartheta(\sigma), \pi)$,}
\end{equation}
and the function $\eta : [0, \pi) \to [0, \infty)$, given by
\begin{equation}\label{PH}
\eta(\varphi)=\frac{ \Psi_{n}(\varphi)}{\sin^{n}\varphi} \quad
\hbox{for $\varphi \in [0, \pi)$,}
\end{equation}
are bijective.
%In particular, if $\vartheta (\sigma) \in
%(0, \pi)$ is defined by
%$$\Psi _n (\vartheta (\sigma) ) = \sigma \Psi _n(\pi),$$
%then the restriction
%$$\xi_\sigma  : (\vartheta (\sigma) , \pi ) \to (0, \infty )$$
%is bijective.
Thus, on defining the (strictly increasing) function $f : [0,\pi
)\to [\vartheta (\sigma),\pi)$ as
\begin{equation}\label{f}
f(\varphi) = \xi_\sigma^{-1}(\eta(\varphi)) \quad \hbox{for $\varphi
\in (0,\pi )$,} \end{equation}
 equation \eqref{conda} is equivalent
to
\begin{equation}\label{explicit}
\vartheta=f(\varphi).
\end{equation}
Note that
\begin{equation}\label{march204}
\varphi < f(\varphi ),
\end{equation}
and hence $(f(\varphi ), \varphi ) \in \Xi $ for $\varphi \in [0,
\pi)$.
%
%
%
%where $f :]0,\pi[\rightarrow]\vartheta (\sigma),\pi[$ is a strictly
%increasing function. This means that the constrained problem:
%\begin{equation}\label{max}
%\sup\{G(\vartheta,\varphi): \vartheta (\sigma)
%<\vartheta<\pi,\>
%0<\varphi<\pi, \text{ (\ref{conda}) holds true}
%\},
%\end{equation}
Altogether, one has that
\begin{equation}\label{maxa}
 \sup_{(\vartheta ,\varphi) \in \Xi}G(\vartheta ,\varphi) = \sup_{\varphi \in [0, \pi)}G(f(\varphi),\varphi).
\end{equation}
The set of critical points of the function  $(0,\pi) \ni \varphi
\mapsto G(f(\varphi),\varphi)$ agrees with the set of the solutions
to the system
\begin{equation}\label{lagrange}
\left\{\substack {\displaystyle\frac{
\Psi_{n-2}(\varphi)}{\sin^{n}\varphi}=\frac{
\Psi_{n-2}(\vartheta)}{\sin^{n}\vartheta}\frac{
(n-1)\cos\vartheta(\Psi_{n-2}(\vartheta)-\sigma
\Psi_{n-2}(\pi))-\sin^{n-1}\vartheta}{(n-1)\cos\vartheta\,\Psi_{n-2}(\vartheta)-\sin^{n-1}\vartheta}
\qquad\qquad\qquad\>\>\\\\
\\
\displaystyle\frac{
(n-1)\Psi_{n-2}(\varphi)-\cos\varphi\,\sin^{n-1}\varphi}{\sin^{n}\varphi}=\frac{
(n-1)(\Psi_{n-2}(\vartheta)-\sigma
\Psi_{n-2}(\pi))-\cos\vartheta\,\sin^{n-1}\vartheta}{\sin^{n}\vartheta},
}\right.
\end{equation}
where $\vartheta $ and $\varphi$ are related as in \eqref{explicit},
and the second equation  is obtained from \eqref{conda}, via
\eqref{ott1}.  On making use of the first equation in
\eqref{lagrange} to rewrite  the second one, and defining the
function $g : (0, \pi ) \to \mathbb R$ as
$$
g (\vartheta)=\frac{\sigma(n-1)\Psi_{n-2}(\pi)}{(n-1)\cos
\vartheta\,\Psi_{n-2}(\vartheta)-\sin^{n-1}\vartheta} \qquad\hbox{for
$\vartheta \in (0, \pi)$,}
$$
system \eqref{lagrange} reads
\begin{equation}\label{lagr}
\qquad\qquad\left\{
\begin{array}{l}
\displaystyle\frac{\Psi_{n-2}(\varphi)}{\Psi_{n-2}(\vartheta)}\frac{\sin^{n}\vartheta}{\sin^{n}\varphi}=1-{g(\vartheta)}{\cos\vartheta}
\\
\\
\displaystyle\frac{\cos\varphi}{\sin\varphi}=\frac{\cos\vartheta}{\sin\vartheta}-\frac{g
(\vartheta)}{\sin\vartheta}.
\end{array}
\right.
\end{equation}
System \eqref{lagr} agrees with  \eqref{sistema}, with $\rho$
replaced by $\sigma$. Also,
\begin{equation}\label{g} g(\vartheta) < 0 \quad \hbox{for
$\vartheta \in (0, \pi ).$}
\end{equation}
Solving  the first equation of \eqref{lagr} for $g(\vartheta)$, and
plugging the resulting expression for $g(\vartheta)$ in the second
equation yield
\begin{equation}\label{necessaria}
G(f(\varphi),\varphi)=\frac1{\cos(f(\varphi)-\varphi)}.
\end{equation}
In conclusion, on defining $F: (0, \pi) \to \mathbb R$ as
$$
F(\varphi)=G(f(\varphi),\varphi)-\frac1{\cos(f(\varphi)-\varphi)}
\quad \hbox{for $\varphi \in (0, \pi)$,}$$
 one has that
\begin{equation}\label{equivalent}
\{\varphi\in (0,\pi) : \varphi \mbox{ is critical for
}G(f(\varphi),\varphi)\}= \{\varphi\in (0,\pi) : F(\varphi)=0\}.
\end{equation}
%
%
%A straightforward computation shows that the function
%$f(\varphi)-\varphi$ has negative derivative with respect to
%$\varphi$ in $]0,\pi[$. Indeed
Computations show that
\begin{equation}\label{march300}
f'(\varphi) -1 = {\frac{n \sin^{n+1}\vartheta
\Psi_n(\varphi)}{\sin^n \varphi [\sin^{n+1} \vartheta  - n \cos
\vartheta (\Psi_n(\vartheta) - \sigma\Psi_n(\pi) )]}}  \left({\frac
{\cos \vartheta}{\sin \vartheta}} -{\frac {\cos \varphi}{\sin
\varphi}} \right)
\end{equation}
for $\varphi \in (0,\pi)$  and $\vartheta = f(\varphi)$. Owing to
the second equation in \eqref{lagr} and to \eqref{g},
$$\left({\frac
{\cos \vartheta}{\sin \vartheta}} -{\frac {\cos \varphi}{\sin
\varphi}} \right)<0.$$
 On the other hand,
$$\sin^{n+1} \vartheta  - n \cos
\vartheta (\Psi_n(\vartheta) - \sigma\Psi_n(\pi) )>0,$$ since the
function on the left-hand side may only achieve a positive minimum.
Hence, equation \eqref{march300} ensures that
\begin{equation}\label{march301}
f'(\varphi) -1 < 0 \qquad \hbox{for $\varphi \in (0,\pi)$.}
\end{equation}
%
%
% {\color{blue}where the inequality is a
%consequence of the second equation in \eqref{lagr} and \eqref{g},
%and of the fact that $\phi(t)=\sin^{n+1} t  - n \cos t (\Psi_n(t) -
%\sigma\Psi_n(\pi) )>0$, for \textcolor[rgb]{0.00,0.59,0.00}{$t \in
%(\vartheta (\sigma ), \pi)$}, in view of the fact that the minimum
%of $\phi(t)$ is achieved at the point $\bar t\in(0,\pi)$ such that
%$\Psi_n(\bar t) = \sigma\Psi_n(\pi)$.}
Owing to \eqref{march204}, one has that
 $0 \leq f(\varphi) - \varphi \leq f(\varphi) \leq \pi$ for
$\varphi \in (0,\pi)$. Therefore,
$$\frac{d}{d\varphi} \bigg(\frac1{\cos(f(\varphi)-\varphi)}\bigg) =
\frac{\sin (f(\varphi) - \varphi )}{\cos^2(f(\varphi) -
\varphi)}(f'(\varphi)-1) <0$$
for $\varphi \in (0, \pi)$. On the
other hand, since $\lim_{\varphi\rightarrow0^+} f(\varphi)=\vartheta
(\sigma)$,
$$\lim_{\varphi\rightarrow0^+}
G(f(\varphi),\varphi)=(n-1)\frac{\Psi_{n-2}(\vartheta
(\sigma))}{\sin ^{n-1}\vartheta (\sigma)},$$ and hence
$$ \lim_{\varphi\rightarrow0^+}
F(\varphi)= \begin{cases} (n-1)\frac{\Psi_{n-2}(\vartheta
(\sigma))}{\sin ^{n-1}\vartheta (\sigma)} - \frac 1{\cos \vartheta
(\sigma)} \quad &
\hbox{if $\sigma \neq \tfrac 12$} \\
-\infty \quad & \hbox{if $\sigma = \tfrac 12$.}
\end{cases}
$$
Next, $\lim_{\varphi\rightarrow\pi^-} f(\varphi)=\pi$. Hence, one
can deduce that $ \lim_{\varphi\rightarrow\pi^-}
G(f(\varphi),\varphi)=(1-\sigma)^{-\frac{n-1}n}$, and therefore
$$
\lim_{\varphi\rightarrow\pi^-} F(\varphi)=\left(\frac1{1-\sigma}\right)^{\frac{n-1}n}-1>0.
$$
In conclusion, if
 $\sigma \in (0, \tfrac 12]$, then  $F$ is a continuously
 differentiable
function on $(0,\pi)$, whose derivative, by \eqref{equivalent}, is
  positive at every point where
$F$ vanishes. Moreover, $\lim_{\varphi\rightarrow0^+} F(\varphi)<0$,
since $\vartheta (\sigma ) \in (0, \tfrac \pi 2]$ and \eqref{g}
holds, and $\lim_{\varphi\rightarrow \pi ^-} F(\varphi)>0$.
\\ If, instead, $\sigma
\in (\tfrac 12 , 1)$, then $\vartheta (\sigma) \in (\pimez , \pi)$,
and since $f (0) = \vartheta (\sigma)$ and $f (\pi ) - \pi =0$,
there exists  a unique $\varphi (\sigma) \in (0, \pi )$ such that $f
(\varphi (\sigma) ) -
 \varphi (\sigma) = \pimez$. Thus, $F$ is a continuously
 differentiable
function on $]0,\pi[\setminus \{\varphi (\sigma)\}$,  whose
derivative is
  positive at every point where
$F$ vanishes, and such that $\lim_{\varphi\rightarrow0^+}
F(\varphi)>0$, $\lim_{\varphi\rightarrow\varphi (\sigma)^-}
F(\varphi)=+\infty$, $\lim_{\varphi\rightarrow\varphi (\sigma)^+}
F(\varphi)=-\infty$, $\lim_{\varphi\rightarrow \pi ^-}
F(\varphi)>0$. \\ In both cases, the set $\{\varphi\in (0,\pi) :
F(\varphi)=0\}$ consists of exactly one point $\varphi _\sigma
 \in (0, \pi)$, whence, by \eqref{equivalent},
\begin{equation}\label{mpoint}
\{\varphi\in (0,\pi) : \varphi \mbox{ is critical for
}G(f(\varphi),\varphi)\}= \{\varphi _\sigma\}.
\end{equation}
One can verify that the derivative of $G(f(\varphi),\varphi)$ is
positive in
  a right neighborhood of $0$, and  negative in a left neighborhood of
  $\pi$. Therefore, the critical
 point $\varphi _\sigma$ for $G(f(\varphi),\varphi)$ is its unique maximum
 point.   On setting $\vartheta_\sigma = f(\varphi _\sigma)$, we
  have thus established \eqref{march202},
where $(\vartheta _\sigma, \varphi _\sigma)$ is the unique solution
in $\Upsilon$  to system \eqref{sistema} with $\rho$ replaced by
$\sigma$. The proof is complete. \qed

\begin{remark} {\rm  In the special case when $n=2$ and $\sigma =1/2$,  the point $\varphi_{1/2}$
satisfies the system
\begin{equation*}
\left\{
\begin{array}{l}
f(\varphi_{1/2})=2\varphi_{1/2} \\
\\
\varphi_{1/2}(1-4\cos^2\varphi_{1/2})=2\sin\varphi_{1/2}\cos\varphi_{1/2}(\frac\pi8-\sin^2\varphi_{1/2}),
\end{array}
\right.
\end{equation*}
where $f$ is defined as in \eqref{f}. In particular, the first
equation tells us that $\varphi_{1/2} = \tfrac 12 \vartheta _{1/2}$,
and hence it entails that the point $P$ in Figure \ref{intersection}
belongs to $\partial \mathbb B^2$.}
\end{remark}

\bigskip
\par\noindent
\textbf{Acknowledgements}.  This research was partly supported by
the research project of MIUR (Italian Ministry of Education,
University and Research) Prin 2012, n. 2012TC7588,  ``Elliptic and
parabolic partial differential equations: geometric aspects, related
inequalities, and applications", and  by GNAMPA of the Italian INdAM
(National Institute of High Mathematics).

\end{document}